\documentclass[12pt]{amsart}

\usepackage[T1]{fontenc}
\usepackage[cp1250]{inputenc}
\usepackage[margin=1in]{geometry}
\usepackage{color}
\usepackage[centertags]{amsmath}
\usepackage{amsfonts}
\usepackage{amsthm}
\usepackage{etex}
\usepackage{newlfont}
\usepackage{amscd}
\usepackage{verbatim}
\usepackage{amsgen}
\usepackage{amssymb}
\usepackage{stmaryrd}
\usepackage{amssymb,amsmath}
\usepackage{amscd}
\usepackage{enumerate}
\usepackage{longtable}
\usepackage{blkarray}
\usepackage{multicol}
\usepackage{graphicx}
\usepackage{colordvi}
\usepackage{verbatim}
\usepackage{cancel}
\usepackage{caption}
\usepackage{epsf}
\usepackage[all,  colour, line, cmtip]{xy}
\usepackage{textcomp}

\usepackage{tikz}
\usepackage{gauss}
\usetikzlibrary{matrix,decorations.pathreplacing,calc}

\newdir{smallhead}{
:(15,1)\dir^{>}
*:(15,-1)\dir_{>}}

\newdir{multihead}{
\dir{smallhead}
*!(0,1.1)\dir{smallhead}
*!(0,-1.1)\dir{smallhead}
*{\POS(-6,0); \POS(0, 1.1) **\crv{~**\dir{-} (-4.5,1.1)}}
*{\POS(-6,0); \POS(0,-1.1) **\crv{~**\dir{-} (-4.5,-1.1)}}
}

\newcommand{\multito}{
\
\begin{xy}
  =(0,0) "SOURCE";
  =(9,0) "TARGET";
  {\ar@{multihead} "SOURCE"; "TARGET"}
\end{xy}
\
}

\newcommand{\multimapsto}{
\
\begin{xy}
  =(0,0) "SOURCE";
  =(9,0) "TARGET";
  {\ar@{|-multihead} "SOURCE"; "TARGET"}
\end{xy}
\
}

\def\kk{k}
\newlength{\defbaselineskip} 
\theoremstyle{plain}
\newtheorem{thm}{Theorem}[section]
\newtheorem{cor}[thm]{Corollary}
\newtheorem{con}[thm]{Conjecture}
\newtheorem{df}[thm]{Definition}
\newtheorem{lema}[thm]{Lemma}
\newtheorem{obs}[thm]{Proposition}
\newtheorem{exer}[thm]{Exercise}
\newtheorem{question}[thm]{Question}

\newtheorem{fact}[thm]{Fact}

\newtheorem{pr}{Algorithm}
\theoremstyle{definition} 

\theoremstyle{definition} 
\newtheorem{ex}{Example}[section] %
\newtheorem{exm}[thm]{Example}
\newtheorem{rem}[thm]{Remark}

 \numberwithin{equation}{section}

\def\p{\mathbb{P}}
\def\r{\mathbb{R}}
\def\z{\mathbb{Z}}
\def\n{\mathbb{N}}
\def\Z{\mathbb{Z}}
\def\N{\mathbb{N}}
\def\c{\mathbb{C}}

\def\q{\mathbb{Q}}
\def\C{\mathbb{C}}

\DeclareMathOperator{\Cut}{Cut}
\DeclareMathOperator{\Supp}{Supp}
\def\fa{\begin{fact}}
\def\kfa{\end{fact}}

\newcommand{\fromto}[2]{#1, \dotsc, #2}

 \DeclareMathOperator{\divi}{div}
\DeclareMathOperator{\Pic}{Pic}
\DeclareMathOperator{\Proj}{Proj}
\DeclareMathOperator{\Imi}{Im}
\newcommand{\set}[1]{\left\{#1\right\}}

\DeclareMathOperator{\Spec}{Spec}
\DeclareMathOperator{\Hom}{Hom}

\def\p{\mathbb{P}}
\def\a{\mathbb{A}}
\def\ob{\begin{obs}}
\def\kob{\end{obs}}
\def\dow{\begin{proof}}
\def\kdow{\end{proof}}

\def\tw{\begin{thm}}
\def\ktw{\end{thm}}
\def\hip{\begin{con}}
\def\khip{\end{con}}
\def\lem{\begin{lema}}
\def\klem{\end{lema}}
\def\ex{\begin{exm}}
\def\prog{\begin{pr}}
\def\kprog{\end{pr}}
\def\wn{\begin{cor}}
\def\kwn{\end{cor}}
\def\uwa{\begin{rem}}
\def\kuwa{\end{rem}}
\def\kex{\end{exm}}
\def\dfi{\begin{df}}
\def\kdfi{\end{df}}
\def\Gf{\mathcal G}
\begin{document}
\title{Selected topics on Toric Varieties}
\author{Mateusz Micha\l ek}
\maketitle
\section{Introduction}
There are now many great texts about or strongly related to toric varieties, both classical or new, compact or detailed - just to mention a few: \cite{Cox, Ful, Oda, Stks, BG}. This is not surprising - toric geometry is a beautiful topic. No matter if you are a student or a professor, working in algebra, geometry or combinatorics, pure or applied mathematics you can always find in it some new theorems, useful methods, astonishing relations.
Still, it is absolutely impossible to compete with the texts above neither in scope, level of exposition nor accuracy. 
We present a review on toric geometry based on ten lectures given at Kyoto University, divided into two parts. The first part is the classical, basic introduction to toric varieties. Our point of view on toric varieties here, is as images of monomial maps. Thus, we relax the normality assumption, but consider varieties as embedded. We hope that this part is completely self-contained - proofs are at least sketched and a motivated reader should be able to reconstruct all details. Further, such an approach should allow the reader not familiar with toric geometry to realize that he encountered toric varieties before. 

The second part deals with slightly more advanced topics. These were chosen very subjectively according to interests of the author. We start by recalling the theory of divisors on toric varieties in Section \ref{sec:div} and their cohomology in Section \ref{sec:Cohom}. Then we present basic results on Gr\"obner degenerations and relations to triangulations of polytopes in Section \ref{sec:GB}, based on beautiful theory by Bernd Sturmfels \cite{Stks}. In these three sections our aim completely changes. We do not present proofs, but focus on methods and examples.

In Section \ref{sec:cutsandsplits} we present toric varieties coming from cuts of a graph. Here, our aim was to prove that a conjecture of Sturmfels and Sullivant implies the famous four color theorem, cf.~Proposition \ref{prop:con=>4col}. This fact is known to experts and the original proof is due to David Speyer - however so far was not published. 

In Section \ref{sec:TVandMat} we present relations of matroids, toric varieties and orbits in Grassmannian. Our focus in on famous White's conjecture and finiteness results related to it.

Mathematical biology provides a source of interesting toric varieties in the case of the so-called phylogenetic statistical models. Such a model can be studied as an algebraic variety by solving the polynomial equations which hold among its marginal probability functions in the appropriate field (e.g. $\r$ or $\C$). An interesting subclass of these models are the group-based models discussed in Section \ref{sec:TVandPhylo}; for each choice of finite graph $\Gamma$ (viewed as the underlying topological structure of the phylogeny) and finite group $G$ there is a well-defined toric ideal $I_{\Gamma, G}$ and affine semigroup $M_{\Gamma, G}$. We describe basic toric constructions and known results in this area.

In section \ref{sec:JB} we very briefly recall the construction of Cox rings and present results of Brown, Buczy{\'n}ski and K{\k{e}}dzierski on their relations to rational maps of varieties.

In the last section \ref{sec:ex} we present some examples related to questions about depth and inner projections of toric varieties.

Throughout the text the reader may find various open problems and conjectures, explicit computations in Macaulay2 \cite{M2}, relying on Normaliz \cite{Normaliz} and 4ti2 \cite{4ti2}. Another great platform for toric computations is Polymake \cite{polymake}.

There are a lot of very interesting, important topics that we do not address and it is impossible even to list all of those. Let us just mention a few (that we regret most to omit): higher complexity  $T$-varieties \cite{altmann2011geometry}, many relations to combinatorics, such as e.g.~binomial edge ideals \cite{ohsugi1999toric}, toric vector bundles \cite{klyachko1990equivariant, payne2009toric}, relations to tropical geometry \cite[Chapter 6]{maclagan2015introduction} and secant varieties \cite{MOZ}.

\section*{Acknowledgements}
I would like to thank very much Takayuki Hibi and Hiraku Nakajima for providing great working conditions at RIMS, Kyoto. I am grateful to Seth Sullivant and David Speyer for sharing the proof of Proposition \ref{prop:con=>4col} (the idea is originally due to David Speyer). I express my gratitude to Jaros{\l}aw Buczy{\'n}ski (section \ref{sec:JB} is based on his advice), Micha{\l} Laso{\'n} (for contributions to section \ref{sec:matOpen}), Christopher Manon (for contributions to section \ref{sec:GBM}), Karol Palka (for pointing out interesting open problems) and Matteo Varbaro (for contributions to section \ref{sec:MatFC}). Further I thank Greg Blekherman, Sijong Kwak, Hyunsuk Moon and Rainer Sinn for discoussions on depth and projections.
\part{Introduction to toric varieties - basic definitions}
\section{The Torus}
In algebraic geometry we study an (affine) algebraic variety $X$, using (locally) rational functions on it - these form a ring $R_X$.
\begin{exm}[Polynomials and Monomials]\label{ex:poly}
Consider the affine space $\C^n$. The associated ring of polynomial functions is $\C[x_1,\dots,x_n]$. We represent a polynomial as:
$$P=\sum_{(a_1,\dots,a_n)\in\N^n}\lambda_{(a_1,\dots,a_n)}x_1^{a_1}\cdots x_n^{a_n},$$
with only finitely many $\lambda_{(a_1,\dots,a_n)}\neq 0$. When $a=(a_1,\dots,a_n)\in \N^n$ we use a multi-index notation $x^a:=x_1^{a_1}\cdots x_n^{a_n}$. Such an expression is called a monomial.
\end{exm}
The main object of these lectures is the complex torus $T=(\C^*)^n$, with the structure of the group given by coordinatewise multiplication. 
On $T$ we have more functions than on the affine space: for $a\in\N^n$ we allow $x^a$ in the denominator.
\dfi[Laurent polynomial]
We define the ring $\C[x_1,x_1^{-1},\dots,x_n,x_n^{-1}]$ consisting of \emph{Laurent polynomials}:
$$P=\sum_{(a_1,\dots,a_n)\in\Z^n}\lambda_{(a_1,\dots,a_n)}x_1^{a_1}\cdots x_n^{a_n},$$
with finitely many $\lambda_{(a_1,\dots,a_n)}\neq 0$. We use the same multi-index notation as in Example \ref{ex:poly}.
\kdfi
This is an example of a more general construction of \emph{localization} that is an algebraic analogue of removing closed, codimension one sets from an affine algebraic variety \cite{AtiMac}. 

Given a map of algebraic varieties $f:X\rightarrow Y$ we obtain an associated map $f^*$ from functions on $Y$ to functions of $X$ (by composition with $Y\rightarrow \C$) \cite{Hart}. Our first aim is to study algebraic maps $(\C^*)^n\rightarrow \C^*$. These correspond to maps of rings $\C[x,x^{-1}]\rightarrow \C[x_1,x_1^{-1},\dots,x_n,x_n^{-1}]$. 
\begin{obs}\label{prop:Tmor}
Every algebraic map $(\C^*)^n\rightarrow \C^*$ is given by $\lambda x^a$ for $\lambda\in\C^*$ and $a\in\Z^n$.
\end{obs}
\dow Fix a lexicographic order on Laurent monomials. For any Laurent polynomial $P$ we denote $LT(P)$ its leading term and $ST(P)$ its smallest term.
Consider the ring morphism associated to the given map.
Suppose $x\rightarrow Q$ and $x^{-1}\rightarrow S$. As $xx^{-1}=1=QS$ we see that  $LT(S)LT(Q)=1=ST(S)ST(Q)$. Hence, $LT(Q)=ST(Q)$. 
\kdow

\wn\label{cor:TorMap} 
Any algebraic map $(\C^*)^n\rightarrow (\C^*)^m$ is given by $(\lambda_1x^{a_1},\dots, \lambda_mx^{a_m})$ for $\lambda_1,\dots,\lambda_m\in \C^*$ and $a_1,\dots,a_m\in\Z^n$. If in addition the map is a group morphism then it is given by Laurent monomials.  
\kwn
\dfi[Characters, Lattice $M$]\label{dfi:M}
Algebraic group homomorphisms $T\rightarrow \C^*$ are called \emph{characters}. They form \emph{a lattice}\footnote{Throughout the text by a \emph{lattice} we mean a finitely generated free abelian group, i.e.~a group isomorphic to some $(\Z^n,+)$.} $M_T\simeq \Z^n$, with the addition induced from the group structure on $\C^*$. Explicitly, given $\chi_1:T\rightarrow \C^*$ and $\chi_2:T\rightarrow \C^*$ we define $$\chi_1+\chi_2:T\ni t\rightarrow \chi_1(t)\chi_2(t)\in \C^*.$$
\kdfi
\ex 
Consider two characters $(\C^*)^2\rightarrow \C^*$:

$$(x,y)\rightarrow x^2y, \quad(x,y)\rightarrow y^{-3}.$$
The first one is identified with $(2,1)\in \Z^2\simeq M_T$ and the second one with $(0,-3)$. Their sum is the character $(x,y)\rightarrow x^2y^{-2}$ corresponding to $(2,-2)\in \Z^2$.
\kex
\dfi[One-parameter subgroups, Lattice $N$]\label{dfi:N}
Algebraic group homomorphisms $\C^*\rightarrow T$ are called \emph{one-parameter subgroups}. They form \emph{a lattice} $N_T\simeq \Z^n$, with the addition induced from the group structure on $T$.
\kdfi
The similarities in Definitions \ref{dfi:M} and \ref{dfi:N} are not accidental. The two lattices are dual, i.e.~$M=\Hom_\Z(N,\Z)$ and $N=\Hom_\Z(M,\Z)$. In other words, there is a natural pairing $M\times N \rightarrow \Z$. Indeed, given a map $m:T\rightarrow \C^*$ and $n:\C^*\rightarrow T$ we may compose them obtaining $\C^*\rightarrow \C^*$ that, by Proposition \ref{prop:Tmor} is represented by an integer.
\begin{exer}
Check that using the identifications $M\simeq \Z^n$ and $N=\simeq \Z^n$ the pairing is given by the usual scalar product. 
\end{exer}
The construction below associates to any monoid $\mathfrak M$ a ring $\C[\mathfrak M]$ known as a monoid algebra. As a vector space over $\C$ the basis of $\C[\mathfrak M]$ is given by the elements of $\mathfrak M$. The multiplication in $\C[\mathfrak M]$ is induced from the monoid action.
\begin{exm}
Consider two monoids: $\mathfrak M_1=(\Z_+^n,+)
$, $\mathfrak M_2=(\Z^n,+)$.
By identifying $x_1^{a_1}\cdots x_n^{a_n}$ with $(a_1,\dots,a_n)$ we see that:
$$\C[\mathfrak M_1]=\C[x_1,\dots,x_n],\qquad \C[\mathfrak M_2]=\C[x_1,x_1^{-1},\dots,x_n,x_n^{-1}].$$
Hence, $\mathfrak M_1$ corresponds to the affine space and $\mathfrak M_2$ to the torus $T$.
More canonically, noticing that elements of the monoid induce functions on the associated variety, we get that the ring of functions on $T$ equals $\C[M_T]$. 
\end{exm}
A map of two tori $f:T_1\rightarrow T_2$ corresponds to a map $\C[M_{T_2}]\rightarrow \C[M_{T_1}]$. By Corollary \ref{cor:TorMap} group morphisms $f$ between two tori correspond to \emph{lattice} maps $\hat f:M_{T_2}\rightarrow M_{T_1}$ or equivalently $N_{T_1}\rightarrow N_{T_2}$.  
\dfi[Saturated sublattice, saturation]\label{def:saturated}
A sublattice $M'\subset M$ is called \emph{saturated} if for every $m\in M$ if $km\in M'$ for some positive integer $k$, then $m\in M'$.

For any sublattice $M'\subset M$ we define its \emph{saturation} by:
$$\{m:\text{there is a positive integer }k\text{ such that }km\in M'\}.$$
\kdfi
\ex 
A sublattice $\{(a,a):a\in\Z\}\subset\Z^2$ is saturated. A sublattice $\{2a:a\in \Z\}\subset \Z$ is not saturated.
\kex
\begin{exer}$ $
\begin{enumerate}
\item Prove that $M'\subset M$ is saturated if and only if there exists a lattice $M_1$ and a lattice map $M\rightarrow M_1$ with kernel $M'$.
\item Show that a saturation of a sublattice is also a sublattice.
\end{enumerate}
\end{exer}
In algebraic geometry, just as we study varieties through function on them, we understand a point of a variety by evaluating functions on it. Hence, to $x\in X$ we associate a map $R_X\rightarrow \C$ from the ring of functions on $X$, that sends $f\rightarrow f(x)$. One of the fundamental theorems of algebraic geometry, Hilbert's Nullstellensatz asserts that we can go the other way round: given a (nonzero) ring morphism $R_X\rightarrow \C$ we can find the (unique) corresponding point $x\in X$.
\ex
Given a point $(a_1,\dots,a_n)\in \C^n$ the corresponding morphism $\C[x_1,\dots,x_n]\rightarrow \C$ sends $x_i\rightarrow a_i$.  

For a point $t\in T$ we have a morphism $f_t:\C[M]\rightarrow \C$. Note that each element of $M$ is not in the kernel (as $f_t(-m)$ is the inverse of $f_t(m)$). Hence, points $t\in T$ correspond to \emph{group} morphisms $M\rightarrow \C^*$. In coordinates, a point $(t_1,\dots,t_n)$ corresponds to a map that assigns $\Z^n\ni(a_1,\dots,a_n)\rightarrow \prod_{i=1}^n t_i^{a_i}\in \C^*$.
\kex
\ob
Given a group morphism of tori $f:T_1\rightarrow T_2$ the image equals a subtorus $T'\subset T_2$. We have a canonical isomorphism $M_{T'}=M_{T_2}/\ker \hat f$.
\kob
\begin{proof}
Consider a subtorus $T'\subset T_2$ with the embedding given by the map $M_{T_2}\rightarrow M_{T_2}/\ker \hat f$. Our aim is to prove that $T'=\Imi f$.

Consider a point of $t\in T_2$ represented by a map $f_t:M_{T_2}\rightarrow \C^*$. We have to show that $f_t$ factors through $\hat f$ if and only if $t\in T'$, i.e.~if and only if $\ker\hat f\subset f_t^{-1}(1)$. The implication $\Rightarrow$ is straightforward.

For the other implication consider the injective morphism $i:M_{T'}=M_{T_2}/\ker \hat f\rightarrow M_{T_1}$. It is enough to show that any morphism $M_{T'}\rightarrow \C^*$ factors through $i$. This follows from the fact that $\C^*$ is a \emph{divisible} group, i.e.~an injective $\Z$ module. More directly, in this case, we can extend the morphism to the saturation of $i(M_{T'})$, one element by one (i.e.~the basis) and then extend to the whole $M_{T'}$. 
\end{proof}
\begin{exer}
Show that the ideal (or even vector space) of equations vanishing on the image of $T_1$ is generated by $\{\chi-1:\chi\in\ker \hat f\}$. 
\end{exer}
Dictionary about torus:
\[
\begin{array}{c|c|c}
\text{geometry}& \text{algebra} &\text{combinatorics}\\
\text{torus }T& \text{algebra }\C[M_T]&\text{lattice }M_T\\
\text{point of }T&\text{surjective ring morphism } \C[M_T]\rightarrow \C&\text{group morphism }M_T\rightarrow\c^*\\
\text{algebraic group morphism}\\ T_1\rightarrow T_2&\text{special ring map }\c[M_{T_2}]\rightarrow\c[M_{T_1}]&\text{lattice map }M_{T_2}\rightarrow M_{T_1}\\
\text{image of such morphism}&\text{kernel of the map}&\text{encoded by kernel}
\end{array}\]

The following theorem is the cornerstone of toric geometry (and representation theory). Recall that a representation of a group $G$ on a vector space $V$ is a group morphism $G\rightarrow GL(V)$. In other words, each element of $G$ provides a linear transformation of $V$ (in a compatible way). In this lectures we additionally assume that $G$ and the morphism $G\rightarrow GL(V)$ are algebraic.
\tw\label{tw:rozklad}
Each representation of a torus $T$ acting on $V$ induces a decomposition $V=\bigoplus_i V_i$, where for each $V_i$ there exists $m\in M$ such that for any $t\in T$ and $v\in V_i$ we have $t(v)=m(t)v$, i.e. $T$ acts on $V_i$ by rescaling the vectors, with different weight for different $i$.
\ktw
\begin{proof}
Consider a map $T\rightarrow GL(V)$. It is represented by:
$$f:t\rightarrow \sum_{m\in M} m(t)A_m,$$
where $A_m$ is just a square matrix of scalars. By $f(t_1t_2)=f(t_1)f(t_2)$ we obtain $A_{m_1}A_{m_2}=0$ for $m_1\neq m_2$ and $A_m^2=A_m$. As $f(1)=id$, we obtain that $x=\sum A_mx$. Hence $V=\bigoplus \Imi A_m$ and on the image of $A_m$ the torus acts by scaling by the character $m$.
\end{proof}
\section{Affine Toric Varieties and Cones}
\dfi[Affine Toric Variety]
An affine toric variety is the closure in $\C^m$ of the map $T\rightarrow (\C^*)^m\subset \C^m$, where the first one is given by group morphism. Equivalently: it is a closure of a subtorus of $(\C^*)^m$ or a closure of an image of a Laurent monomial map. In particular, the affine toric variety can be identified with a set of $m$ points $S\subset M_T$.
\kdfi
\ex \label{ex:cusp}
Consider a map $\C^*\rightarrow \C^2$ given by $t\rightarrow (t^2,t^3)$. The associated toric variety is a \emph{cusp} -- the zero locus of $x^3-y^2$. The two points representing this toric variety are $\{2,3\}\subset \Z$.
\kex
\uwa
The definition of (affine) toric variety is often different in other sources and requires the variety to be \emph{normal}.
\kuwa
\tw\label{tw:binomials}
Let $\tilde S$ be the monoid generated by $S$ in $M$. The toric variety $X$ associated to $S$ is isomorphic to $\Spec \C[\tilde S]$. The ideal of $X\subset \C^m$ is linearly spanned by such binomials $y_1^{b_1}\dots y_m^{b_m}-y_1^{c_1}\dots y_m^{c_m}$ that $\sum b_i s_i=\sum c_i s_i$ in $M$, for $b_i,c_i\in\n$.
\ktw
\dow 
It is enough to provide the stated description of the ideal.

It is straightforward that binomials of the given form belong to the ideal. 

Consider any $f(y)$ in the ideal.
We will prove that $f(y)$ is a linear combination of binomials of the given form, inductively on the number of monomials appearing (with nonzero coefficients) in $f$. If there are no monomials (i.e. $f=0$), the statement is obvious.

Choose a monomial $m$ appearing in $f$. If we substitute $y_i\rightarrow x^{a_i}$, where $a_i$ represent characters from $S$, we know that $f$ is zero, as it vanishes on the image. In particular, after the substitution for $m$ (that will remain a monomial, but now in $x$), the obtained monomial must cancel with some other monomial. This other monomial must come from $m'$ appearing in $f$. But the fact that after the substitution they cancel, is equivalent to $m-m'$ being the binomial of the given form. Hence, we may subtract this binomial (with appropriate coefficient), reducing the number of monomials appearing in $f$.
\end{proof}
\uwa
The proof of the Theorem \ref{tw:binomials} does not depend on the field. In fact, the binomials that generate the ideal of $X$ do not depend on the field, cf.~Lemma \ref{lem:fieldindependent}.
\kuwa

\tw\label{tw:afftoric}
An affine variety on which the torus $T$ acts and has a dense orbit is an affine toric variety.
\ktw
\dow 
The algebra $R_X$ of the variety embeds into $\C[M]$ as the morphism $t\rightarrow tx$ is dominant for (general) $x\in X$. We claim that $R_X$ is linearly spanned by elements of $M$. Indeed, consider any $g\in R_X$. The torus acts on $R_X$, in particular on $g$. Consider the (finite dimensional) vector space spanned by all $Tg$. By Theorem \ref{tw:rozklad} all characters $\chi_i$ appearing in $g=\sum_i c_i\chi_i$ (with nonzero coefficients) must belong to $R_X$.

As $R_X$ is finitely generated, the monoid of characters $\tilde S$ in $R_X$ is finitely generated, with generators providing the embedding in the affine space.

It is worth noticing that the dense torus orbit must be in fact also a torus.
\kdow
The following definition extends Definition \ref{def:saturated} to monoids.
\dfi[Saturated monoid]\label{def:SatMon}
A monoid $\tilde S\subset M$ is saturated (in $M$) if and only if $km\in\tilde S$ for some $k\in \N_+$, $m\in M$ implies $m\in \tilde S$.
\kdfi 
We now introduce a quite subtle notion of \emph{normality}. 
\dfi[Integrally closed, Normal]\label{def:intcl+norm}
We say that a ring $A\subset B$ is integrally closed in $B$ if for any monic (i.e.~with the leading coefficient equal to $1$) polynomial $f\in A[x]$ if for some $b\in B$ we have $f(b)=0$ then $b\in A$.   

We say that an integral ring $A$ is normal if it is integrally closed in its ring of fractions.
\kdfi
At this point the definition of normality may look artificial. It turns out that when a ring $R_X$ is normal then the variety $X$ is not 'too singular'. In particular, if $X$ is smooth then $R_X$ is normal. Normality turns out to play a crucial role in many branches of mathematics. Below we will see it appears naturally in toric geometry. Later, we will point out connections to the properties of divisors and finally we will show relations e.g.~to matroid theory. 

\dfi[Cone]\label{def:cone}
By a cone $C$ in a lattice $M$ (resp.~vector space $V$ over $\r$ or $\q$) we mean a subset containing $0$ and closed under any nonnegative linear combinations:

If $\sum \lambda_i c_i\in M$ for $\lambda_i\in \r_+$, $c_i\in C$ then $\sum \lambda_i c_i\in C$.

A cone is called \emph{polyhedral} if it is finitely generated (using nonnegative linear combinations) and
\emph{rational} if its generators are lattice points. 
\kdfi
Note that positive even integers do not form a cone in $\Z$ but they do form a cone in $2\Z$.
\tw\label{tw:charnormal}
The affine toric variety $X$ is normal if and only if the associated monoid $\tilde S$ is saturated in the lattice that it spans.

A saturated monoid is a cone. Every finitely generated cone is finitely generated as a monoid.
\ktw
\dow 
First let us prove that if $X$ is normal then $\tilde S$ is saturated. Consider any point $kc\in \tilde S$. We want to prove that $c\in \tilde S$. Let $M$ be the lattice spanned by $\tilde S$.
To improve notation, for $m\in M$ let $\chi_m$ be a corresponding character. Consider a polynomial $f(X)=X^k-\chi_{kc}$ with coefficients in the algebra of $X$. Due to the normality of $X$ we know that $\chi_c$ is also in the algebra. Hence $c\in \tilde S$.

It remains to prove that if $\tilde S$ is saturated, then $\c[\tilde S]$ is normal. First note that the quotient field of $\c[\tilde S]$ is equal to the quotient field of $\c[M]$. As the torus is smooth, its algebra is normal. (One can also prove it by noticing that its algebra is a UFD  - as it is a localization of the polynomial ring.) Consider any monic polynomial $f\in\c[\tilde S][x]$. Suppose that $g$ is in the quotient field and satisfies the equation $f(g)=0$. From the normality of $\c[M]$ we know that $g\in \c[M]$. We can act on the equation $f(g)$ by any point $t$ of the torus $T$. The action of $t$ on $f$ gives a monic polynomial with coefficients in $\c[\tilde C]$. Hence the action of $T$ on $g$ gives polynomials that are in the normalization of $\c[\tilde S]$. Considering the action of $T$ on the space of such polynomials, by Theorem \ref{tw:rozklad}, we conclude that all the characters-monomials appearing in $g$ with nonzero coefficient must be in the normalization of $\c[\tilde S]$.
Thus we can assume that $g\in M$.
Suppose that $f$ is of degree $d$. Notice that $f(g)=0$ implies that $dg=d'g+c_0$ for some integer $0\leq d'<d$ and $c_0\in \tilde S$, as the character $\chi_{dg}$ must reduce with some other character. Thus $(d-d')g\in \tilde S$ and by normality $g\in \tilde S$.

For the cone: generators of the monoid belong to the $\{\sum \lambda_i v_i: 0\leq\lambda_i\leq 1\}$, where $v_i$ are generators of the cone.
\kdow
Hence, we have an equivalence of affine normal toric varieties with the torus $\c[M]$ action and (finitely generated, full dimensional) cones in $M$.

\ex 
If we consider all lattice points $(a,b)\in \N^2$ such that $a\leq b\sqrt 2$ then we obtain a monoid that is not finitely generated.
\kex
\ex \label{ex:afftor}
\begin{enumerate}\
\item Looking at Example \ref{ex:cusp} we see that the monoid spanned by $\{2,3\}$:
\begin{itemize}
\item as a group spans $\z$,
\item is not saturated as $2\cdot 1=2$.
\end{itemize}
Indeed, a normal variety is smooth in codimension one (i.e.~the singular locus must have codimension at least two). Our variety is a curve, hence is normal if and only if it is smooth. Note that $0$ is the singular point of the cusp.
\item We already know that the positive orthant gives rise to the affine space and whole lattice $M$ to the torus.
\item Consider the cone generated (as a cone) by $(1,0),(1,2)$. There is one more monoid generator: $(1,1)$. Hence, our variety is realized as (the closure of) the monomial map:
$$(t_1,t_2)\rightarrow (t_1,t_1t_2,t_1t_2^2).$$
We have an integral linear relation: $(1,0)+(1,2)=2\cdot (1,1)$ (all other are generated by this one). Hence, the ring of our variety is: $\C[x,y,z]/(xz-y^2)$.
\end{enumerate}
\kex
\begin{exer}
Prove that each finitely generated cone has a unique, finite minimal set of generators.

Hint: consider all elements that do not have a (nontrivial) presentation $c=c_1+c_2$. 
\end{exer}
Note that we do not have to choose the minimal set of generators of the cone. Continuing Example \ref{ex:afftor} point iii) we may consider a fourth generator, e.g.~$(2,3)=(1,2)+(1,1)$. This provides an embedding in a four dimensional affine space and an isomorphic ring with a different presentation:
$$\C[x,y,z,t]/(xz-y^2,t-yz).$$
\ex 
\begin{enumerate}\
\item The map given by all monomials of degree $r$ (in $n$ variables) is called the r-th Veronese. 
The associated monoid consists of all points in the positive quadrant with sum divisible by $r$. Question: Is this toric variety normal?
\item Consider $k$ groups of (distinct) variables, the $i$-th group consisting of $a_i$ variables. The map given by all monomials (of degree $k$) that are of degree one with respect to each group is called the Segre map.
\end{enumerate}
Both of the above examples are usually considered in projective setting - we will be coming back to them.
\kex
\begin{exer}
Consider the map $(x,y)\rightarrow (x^2,y^2)$. What is the associated toric variety? Is it normal?
\end{exer}
Dictionary about affine toric varieties:
\[
\begin{array}{c|c|c}
\text{geometry}& \text{algebra} &\text{combinatorics}\\
\text{affine toric variety}& \text{prime binomial ideal}&\text{f.g.~submonoid in a lattice}\\
\text{point}&\text{surjective morphism to }\C&\text{semigroup map }S\rightarrow (\C,\cdot)\\
\text{normal}&\text{normal ring}&\text{cone}\\
\text{toric embedding}&\text{special monomial map}&\text{choice of generators of the monoid}\\
\text{hypersurfaces that contain}&\text{defining ideal}&\text{integral relations among generators}
\end{array}
\]
\section{Projective Toric Varieties}
\dfi[Projective Toric Variety]
A projective toric variety is the closure in $\p^m$ of the map $T\rightarrow (\C^*)^m\subset \p^m$, where the first one is given by a group morphism and the inclusion can be regarded as the locus of points with nonzero coordinates.

Equivalently: it is a closure of a subtorus of $(\C^*)^m$ or a closure of an image of a Laurent monomial map in a projective space.
\kdfi
\uwa
As the points in $\p^m$ are regarded up to scalar as $m+1$-tuples of complex numbers, the projective toric variety can be described as a closure of the map given by $m+1$ characters. 
\kuwa

Note that if $m_1,\dots,m_{m+1}$ are monomials parameterizing the projective toric variety $X$, it is easy, using an additional variable $x_0$, to parameterize the affine cone over $X$:
$$(x_0,\dots,x_n)\rightarrow(x_0m_1,\dots,x_0m_{m+1})\in\a^{m+1}.$$
Hence, it is 'better' to represent $m_i$ not in the lattice $\Z^n$, but $\Z^{n+1}$, by the inclusion $\Z^n\ni m\rightarrow (m,1)\in \Z^{n+1}$.
\dfi[Normal polytope]
A lattice polytope (i.e.~a polytope whose vertices are lattice points) $P\subset M_\r$ is called \emph{normal} (in the lattice $M$) if and only if for any $k\in \n$ all lattice points in $kP$ are sums of $k$ (not necessary distinct) lattice points of $P$.
\kdfi
\dfi[Projective normality]
A projective algebraic variety is called projectively normal if and only if the affine cone over it is normal.
\kdfi
Thus, by Theorem \ref{tw:charnormal} we obtain the following.
\tw\label{tw:normalpoly}
A projective toric variety is projectively normal if and only if the associated parameterizing monomials $S\in \Z^{n+1}=<S>$ (at height one) generate (as a monoid) all lattice points in the cone generated by them. 

\ktw
In the situation of the theorem above the lattice points must be all lattice points of a convex polytope, but this is not enough! 
\begin{exer}
Show that a projective toric variety is projectively normal if and only if the associated polytope is normal.
\end{exer}

A set of points (usually integral points in a polytope) $S\subset \Z^{m}$ define a projective toric variety, with a projective embedding in a projective space with a distinguished torus $T'$. The dense torus orbit is the intersection of $T'$ with the variety. Notice that coordinates of the ambient projective space correspond to the points in $S$. Hence, to choose an affine open chart we have to choose one point and 'set it equal to one'. This corresponds to dividing parameterizing monomials by the chosen one and division corresponds to shifting (subtracting the chosen point). In such a way we obtain the open affine toric variety. 

\dfi[Very ample polytope]
A lattice polytope $P\subset M_\r$ is called \emph{very ample} (in the lattice $M$) if and only if all lattice points in $kP$ are sums of $k$ (not necessary distinct) lattice points of $P$ for $k$ large enough.
\kdfi
\tw\label{tw:very ample}
A set of lattice points in a polytope defines a normal projective toric variety if and only if the polytope is very ample (in the lattice it spans).


\ktw
The proof follows from the exercise below.
\begin{exer}
Suppose that $P$ is a lattice polytope that spans a lattice $M$. Show that the following are equivalent:
\begin{itemize}
\item $P$ is very ample,
\item for any lattice point $m\in P$ the monoid spanned by $P-m$ is saturated,
\item for any vertex $m$ of $P$ the monoid spanned by $P-m$ is saturated.
\end{itemize}
\end{exer}
\begin{exer}\
\begin{itemize}
\item Show that the unit $r$-dimensional simplex corresponds to $\p^r$. 
\item The map $\p^r\rightarrow \p^{{d+r}\choose{r}-1}$ given by all monomials of degree $d$ is called the \emph{Veronese} embedding. What is the corresponding polytope?
\item The map $\p(V_1)\times\dots\times\p(V_n)\rightarrow \p(V_1\otimes\dots\otimes V_n)$ defined by $$[v_1]\times\dots\times[v_n]\rightarrow [v_1\otimes\dots\otimes v_n]$$ is called the \emph{Segre} embedding. What is the corresponding polytope?
\item What are the defining equations in the two examples above? Are the varieties projectively normal?
\end{itemize}
\end{exer}
\ex\
\begin{itemize}
\item The set $\{0,1,3,4\}\subset\Z$ defines a smooth projective toric variety that is not projectively normal.
\item The $3$-dimensional polytope with vertices $$(0,0,0),(0,0,-1),(0,1,0),(0,1,-1),(1,0,0),(1,0,-1),(1,1,3),(1,1,4)$$ is very ample, but not normal.
\end{itemize}
\kex
\dfi[Smooth polytope]
An $n$-dimensional lattice polytope $P$ is  \emph{smooth} if 
\begin{itemize}
\item for any vertex $v$, there are exactly $n$ edges adjacent to $v$ with lattice points nearest to $v$ denoted by $v_1,\dots,v_n$ and
\item the vectors $v-v_1,\dots,v-v_n$ form a basis of the lattice spanned by $P$. 
\end{itemize}
\kdfi
\uwa
A lattice polytope $P$ defines a smooth projective toric variety if and only if $P$ is smooth.
\kuwa
Below we present a famous, central conjecture in toric geometry. The first part is due to Oda and the second to Bogvad. 
\hip
A smooth polytope is normal. The associated toric variety is defined by quadrics.
\khip
\ob
The degree of a projective toric variety $X$ represented by a very ample polytope $P$ spanning a lattice $M$ equals the (normalized) volume of $P$.
\kob
\uwa
The function $n\rightarrow |nP\cap M|$ is a polynomial known as the Ehrhart polynomial. 
\kuwa
\dow 
By definition, the degree of $X$ is (up to $\frac{1}{(\dim P)!}$) the coefficient of the leading term of the Hilbert polynomial (that to $n$ associates the dimension of degree $n$ part of $R_X$). However, this space is spanned by lattice points in $nP$. 

To estimate the number of lattice points in $nP$ and relate it to a volume we cover/inscribe in $P$ small cubes $C_i$. The statement relating the volume of a cube $C_i$ with the estimate of the degree of the function $n\rightarrow |nC_i\cap M|$ easily reduces to one dimension and is an easy exercise. 
\kdow

Let us return to an affine covering of the toric variety defined by $P$. 
\begin{exer}
Suppose a set of characters $S$ defines a projective toric variety $X_S$. Show that, under this embedding, $X_S$ is covered by principal affine subsets corresponding to vertices of the convex hull (in $M_\r$) of $S$.
\end{exer}
The previous exercise shows that the monoids generated by $S-v$ are of particular interest. When $S$ is formed by lattice points of a polytope these monoids are cones precisely when the polytope is very ample. 
\subsection{General Toric Varieties}
We briefly describe general toric varieties, that do not have to be affine or projective. This is a classical topic, central in toric geometry, well-described in many books and articles. However, in this review, general toric varieties are a side topic - we focus our interest on projective ones.
\dfi[Fan]
A (finite) collection of (rational, polyhedral) cones in (a vector space over) a lattice $N$ is called a fan if it is closed under taking faces and intersections.
\kdfi
\dfi[Dual cone]
Fix a cone $\sigma\subset M$ (resp.~$M_\r$). We define the dual cone $\sigma^\vee\subset N=M^*$ (resp.~$N_\r$) consisting of those elements $n$ such that for any $m\in \sigma$ we have $(n,m)\geq 0$.  
\kdfi
If we look at the family of cones $P-v$ (where $v$ runs over vertices of $P$) in $M$ we do not see any structure.
\begin{exer}
Show that $(P-v)^\vee$ (together with their faces) form a fan covering whole $N$. 
\end{exer}
Fans covering $N$ are called \emph{complete}. Not every complete fan comes from a polytope $P$, as we will see in Section \ref{sec:div}. However, there is a general construction that to a fan $\Sigma$ in $N$ associates a normal toric variety $X(\sigma)$. The idea is to glue together the affine toric varieties $\Spec\C[\sigma^\vee]$ for all $\sigma\in \Sigma$. Precisely, consider the cones $\sigma_3=\sigma_1\cap\sigma_2$. We have induced open embeddings $\Spec\C[\sigma_3^\vee]\rightarrow\Spec\C[\sigma_i^\vee]$ for $i=1,2$ that allow the gluing. This gives a normal algebraic variety with an action of the torus $\Spec\C[N^\vee]$ with a dense orbit. 

Conversely, given a normal algebraic variety with an action of a torus $T$ and a dense orbit one can prove that it is represented by a fan. The easy part of the proof is to consider torus invariant affine open subvarieties and show that they glue in a good way. However, the nontrivial part of the theorem relies on the fact that torus invariant affine subvarieties cover the whole variety. This is a theorem of Sumihiro \cite{sumihiro1974equivariant, kempt339toroidal}.
\section{Affine Toric Varieties - combinatorics and algebra}
\subsection{Orbit-Cone correspondence}\label{sec:orbit-cone}
For an affine toric variety corresponding to a cone $C$ the faces of $C$ correspond to orbits of the torus acting on it. Let us present this correspondence in details. We fix a finitely generated monoid $C$ in a lattice $M$ and its generators $\chi_1,\dots,\chi_k\in C$. We know that:
\begin{itemize}
\item the dense torus orbit of $X$ contains precisely those points that have all coordinates different from zero,
\item the character lattice of the torus acting on $X$ is equal to the sublattice of $M$ spanned by $C$.
\end{itemize}
We will generalize this to other orbits.
\tw\label{tw:orbitcone}
Assume that $C$ is a cone. Each orbit will be indexed by a face $F$ of the cone. The face $F$ distinguishes a subset $I_F$ of indices from $\{1,\dots,k\}$ such that $i\in I$ if and only if $\chi_i\in F$. The orbit corresponding to $F$ can be characterized as follows:
\begin{enumerate}[  1)]
\item the orbit contains precisely those points that have got coordinates corresponding to $i\in I_F$ different from zero and all other equal to zero,
\item the orbit is a torus with a character lattice spanned by elements of $F$,
\item the closure of the orbit is a toric variety given by the cone $F$,
\item each point of the orbit is a projection of the dense torus orbit onto the subspace spanned by basis elements indexed by indices from $I_F$,
\item the inclusion of the orbit in the variety is given by a morphism of algebras $\c[C]\rightarrow\c[F]$. This morphism is an identity on $F\subset \c[C]$ and zero on $C\setminus F$.
\end{enumerate}
Note that each orbit will contain a unique distinguished point given by the projection of the point $(1,\dots,1)\in \c^k$.
\ktw
\dow
As in case of the torus we can identify the points of $X$ with monoid morphisms $C\rightarrow(\c,\cdot)$. 

Point 1): Fix any point $x\in X$. Note that for any $c_1,c_2\in C$ if $(c_1+c_2)(x)\neq 0$ then $c_1(x),c_2(x)\neq 0$. Hence, the characters $\chi\in C$ such that $\chi(x)\neq 0$ must form a face $F$ of $C$. Thus $x$ distinguishes a subset of indices $I_F\subset\{1,\dots,k\}$. Of course the set of points with nonzero coordinates indexed by $I_F$ and other coordinates equal to zero in $X$ is invariant with respect to the action of the torus acting on $X$. So to prove 1) it is enough to prove that all these points are in one orbit. The point $x$ represents a morphism $C\rightarrow(\c,\cdot)$ that is nonzero on $F$ and zero on $C\setminus F$. Consider the restriction of this morphism to $F$. As it is nonzero it can be extended to a morphism $ M'\rightarrow \c^*$, where $M'$ is a sublattice generated by $F$. Next, as $M'$ is saturated, we can extend this morphism to the lattice $M$ generated by $C$. Thus we obtain a morphism $f:M\rightarrow \c^*$ that agrees with the one representing $x$ on $F$. Note that $f$ represents a point $p$ in the dense torus orbit of $X$. By the action of $p^{-1}$ on $x$ we obtain a point given by a morphism that associates one to elements from $F$ and zero to elements from $C\setminus F$. Thus we have proved 1). Moreover, we showed that each orbit contains the distinguished point. 

Point 2) follows, as morphisms that are nonzero on $F$ and zero on $C\setminus F$ are identified with morphisms from $M'$ to $\c^*$. 

Point 3): We already know that the orbit is a torus with the lattice generated by $F$. This torus is the image of the torus $\Spec \c[M]$ in $\c^k$ by characters from $I_F$ and all other coordinates equal to zero. We see that the orbit corresponding to $F$ is contained in the affine space spanned by basis elements indexed by indices in $I_F$. In fact, by construction it is the image of $\Spec \c[M]$ by characters $\chi_i$, such that $i\in I_F$. The closure of this torus is exactly given by $\Spec \c[F]$, as generators of the monoid $C$ contained in $F$ are generators of $F$. 

Point 4) is obvious, as the point $p$ constructed in the first part of the proof projects to $x$.

Point 5) is a consequence of the other points.
\kdow
Looking at a cone over a polytope we see that orbits correspond to faces of polytope (the whole polytope to the dense torus orbits, facets to torus invariant Weil divisors,\dots, vertices to torus invariant points - the containment of closures of torus orbits is as you can see on the polytope).
\subsection{Toric ideals}
We now characterize ideals $I\subset\c[x_1,\dots,x_n]$ that define toric varieties under a toric embedding. 
\tw\label{tw:bionomial prime}
The ideal $I$ defines an affine variety that is a closure of a subtorus in $\C^n$ if and only if $I$ is a binomial (i.e.~generated by binomials) prime ideal. 
\ktw
\dow 
In Theorem \ref{tw:binomials} we proved the forward implication.

Suppose $I$ is prime and binomial. As $(1,\dots,1)\in V(I)$ we see that the intersection $X_0:=V(I)\cap (\c^*)^n$ is a nonempty variety and hence dense in $V(I)$. Note that a restriction of a binomial to $(\c^*)^n$ is (up to invertible element) of the form $m-1$, where $m$ is a character. As $X_0$ is reduced and irreducible the set $M':=\{m: m-1\in I(X_0)\}$ is a saturated sublattice of a lattice $M_{(\c^*)^n}$. The map $\c[x_1^{\pm 1},\dots,x_n^{\pm 1}]/I(X_0)\rightarrow \c[ M_{(\c^*)^n}/M']$ is an isomorphism. Hence, $X_0$ is a subtorus and $V(I)$ and its closure is a toric variety.
\kdow
\uwa
In general, it can be very hard to say if $I$ defines a toric variety (with a different embedding). In particular, there exists a smooth hypersurface $H$ in $\c^5$ given by $x+x^2y+z^2+t^3$. It is a product of $\C$ and a so-called the Koras-Russel cubic. It is an open problem if it is isomorphic with $\c^4$. Further, it is not known if there exists an automorphism of $\c^5$ that would linearize $H$.  
\kuwa
\part{Topics on Toric Varieties}
\section{Divisors}\label{sec:div}
\subsection{Weil Divisors}
Given any variety $X$ one may consider the set of all irreducible codimension one subvarieties and the free abelian group $D(X)$ generated by it. An integral combination of such subvarieties is called a \emph{Weil divisor}. This group is far too large, thus one regards it modulo an equivalence relation. From now on assume $X$ is regular in codimension one, i.e.~the codimension of the singular locus is at least two (e.g.~$X$ is normal). In such a case to a rational function $f$ on $X$ we can associate a divisor $\divi(f)\in D(X)$ - see e.g. \cite[Chapter 6]{Hart}.
\dfi[Class group]
The class group $Cl(X)$ is defined as $D(X)$ modulo the subgroup generated by all $\divi(f)$. 
\kdfi
Let $X$ be a normal toric variety defined by a fan $\Sigma$ (we may think $\Sigma$ is a normal fan of a polytope - for a general construction we refer to \cite{Cox}). The beautiful thing is that to compute $Cl(X)$ we do not have to consider all divisors (just those invariant by the torus action) and not all rational functions (just characters of the torus).

Precisely we recall that rays in $\Sigma$ correspond to codimension $1$ orbits of the torus (hence their closures are codimension one torus invariant subvarieties). For a ray $u\in \Sigma_1$ we abuse the notation and denote also by $u$ the first nonzero lattice point on that ray. We note that for $m\in M$ we have a pairing $\langle m,u\rangle$ that is equal to zero if and only if $m$ does not vanishes or have a pole on the variety $D_u$ associated to $u$ (think e.g.~in terms of the coordinates as in Chapter \ref{sec:orbit-cone}). Indeed, $\divi(m)=\sum_{u\in\Sigma_1}\langle m,u\rangle D_u.$
\tw
We have an exact sequence:
$$M\rightarrow \bigoplus_{u\in\Sigma_1} \z D_u\rightarrow Cl(X)\rightarrow 0.$$
\ktw
\begin{exm}\label{ex:Weil}\
\begin{enumerate}
\item As $\p^n$ corresponds to a simplex the rays of the normal fan are $e_1,\dots,e_n,-e_1-\dots-e_n\in N$. The exact sequence above becomes:
$$0\rightarrow\Z^n\rightarrow\Z^{n+1}\rightarrow \Z\rightarrow 0.$$
\item As $\p^1\times\p^1$ corresponds to a square the rays of the normal fan are $e_1,-e_1,e_2,-e_2$. The exact sequence becomes:
$$0\rightarrow\Z^2\rightarrow\Z^4\rightarrow\Z^2\rightarrow 0.$$
\item (Rational normal scroll, Hirzebruch surface) Consider a polytope that is the convex hull of $0,f_1,f_2,rf_1+f_2$. It has four facets corresponding to four rays in the dual fan $e_1,e_2,-e_2, -e_1+re_2$. We obtain $Cl(X)=\z^2$. 
\item Consider the cone generated by $f_1,f_1+f_2,f_1+2f_2$. The dual cone has ray generators $e_2,2e_1-e_2$. We obtain:
$$0\rightarrow\Z^2\rightarrow\Z^2\rightarrow\Z/2\Z\rightarrow 0.$$
\end{enumerate}
\end{exm}
\subsection{Cartier divisors}
A \emph{Cartier} divisor is represented by a(n open) covering $U_i$ of $X$ and rational functions $f_i$ such that $f_i/f_j$ is well defined and nonzero on $U_i\cap U_j$. More formally it is a global section of the quotient sheaf of (the sheaf of) invertible rational functions modulo (the sheaf of) invertible regular functions.
Again we consider them modulo divisors given by one rational function (and the covering consisting just of $X$). The quotient is the \emph{Picard group} $\Pic(X)$ of $X$.

As before, for toric varieties we only have to consider $U_i$ affine toric (hence represented by cones) and rational functions $f_i$ given by characters of the torus.
\tw
Let $\Sigma$ be a fan in which each cone is contained in a full dimensional cone. Consider $T$-invariant Cartier divisors $CDiv_T(X)$ given by:
\begin{itemize}
\item A continuous, locally linear function represented by $m_\sigma\in M$ for each maximal dimensional cone $\sigma\in\Sigma$.
\end{itemize}
We have an exact sequence:
$$M\rightarrow CDiv_T(X)\rightarrow\Pic(X)\rightarrow 0.$$
(Without the assumption that maximal cones are full dimensional the exact sequence also holds, but we must define $CDiv_T$ taking into account that $m_\sigma$ is a class in $M/(\sigma^\perp)$.)
\ktw
\begin{exer}
Compute the Picard group in cases 1)-3) from Example \ref{ex:Weil}. Can you see a general theorem here?

What happens in case 4)?
\end{exer}
In general on a normal variety there is an injection $\Pic(X)\rightarrow Cl(X)$. In toric case, we just evaluate the continuous linear function on the ray generators.

Let us now determine the global sections of a line bundle associated to a Cartier divisor $D$ represented by a piecewise linear function $f$ on a fan $\Sigma$. As everything is torus equivariant, the sections will be spanned by characters. Let us fix a character $m\in M$. We ask when $\chi_m\in \mathcal{O}(D)$. Let us fix a maximal cone $\sigma\in \Sigma$. Here $D$ is represented by a rational function $m_\sigma\in M$. Hence, we ask when the product (that is the sum in the lattice) $m+m_\sigma$ is well-defined on $\Spec \C[\sigma^\vee]$. This is clearly if and only if $m+m_\sigma\in \sigma^\vee$. 

Caution: there are different conventions about signs of $f$, $m_\sigma$ etc.

\tw
The global sections of $\mathcal{O}(D)$ are spanned by such $m$ that $m+f$ is nonnegative of $|\Sigma|$. Given a Weil divisor $\sum_{u\in\Sigma_1} a_uD_u$ we define a polyhedron:
$$P_D:=\{m\in M_\q: \langle m,u\rangle \geq -a_u\text{ for all }u\in\Sigma_1\}$$
\ktw
From now on suppose $\Sigma$ is a fan with all maximal cones full dimensional. 
\dfi[Convex, strictly convex]

We say that a piecewise linear function $h$ represented by $m_\sigma$ on cone $\sigma\in\Sigma$ is \emph{convex} if each $m_\sigma\leq h$ (as functions on $|\Sigma|$). We say it is \emph{strictly convex} if the inequality is strict outside $\sigma$. 
\tw
A Cartier divisor $D$ represented by a function $f$ is:
\begin{enumerate}
\item globally generated if and only if $h$ is convex,
\item ample if and only if $h$ is strictly convex,
\item very ample if and only if $h$ is strictly convex and $P_D$ is very ample.
\end{enumerate}
\ktw
\dow 
1) and 3) follow by local affine analysis.
2) is a consequence of 3).  
\kdow
\kdfi

\section{Cohomology}\label{sec:Cohom}
We start with general remarks on cohomology of
line bundles on toric varieties.
Let $\Sigma$ be a fan in $N_\r = \r^n$ with rays $r_1, \dots, r_s$. 
For a Weil divisor $D=\sum_{i=1}^s a_i D_i$ and $m\in M$ we define $$C_{D,m}:=\cup_{\sigma\in \Sigma} Conv({r_i}:r_i\in\sigma, \langle m,r_i\rangle <-a_i).$$ For a Cartier divisor $D_h$ represented by a (piece-wise linear) function $h$ and $m\in M$ we define
$$C_{D_h,m}:=\{u\in |\Sigma|:\langle m,u\rangle <h(u)\}.$$

\ob[\cite{Cox} Section 9.1]\label{prop:coh1}
Note that for a Weil or Cartier torus invariant divisor $D$, each cohomology is $M$ graded. We have:
$$H^p(X_\Sigma, \mathcal{O}(D))_m\simeq \tilde H^{p-1}(C_{D,m})\text{ for a Weil divisor }D,$$
$$H^p(X_\Sigma,D_h)_m=\tilde H^{p-1}(C_{D_h,m})\text{ for a Cartier divisor }D_h,$$
where $\tilde H$ is the reduced cohomology.
\kob

There is an 'Alexander dual' method of computing the cohomology. Assume $\Sigma$ is complete and smooth.
For $I\subset\{1, \dots ,m\}$ let $C_I$ be the simplicial complex generated by sets $J\subset I$
such that $\{r_i : i \in J \}$ form a cone in $\Sigma$. For $a= (a_i : i = 1, \dots ,m)$ let us define
$\Supp(a) := C_{\{i:a_i\geq 0\}}$.
\ob[\cite{borisov2009conjecture}]\label{prop:coh2}
The cohomology $H^j (X_\Sigma, L)$ is isomorphic to the direct sum over all $a =
(a_i : i = 1, \dots,s)$ such that $\mathcal{O}(\sum_{i=1}^s
a_i D_{r_i} )=
L$ of the $(n-1-j )$-th reduced homology of the
simplicial complex $\Supp(a)$.
\kob

To compare Proposition \ref{prop:coh1} and \ref{prop:coh2}, let us first look at $H^0(X_\Sigma,\mathcal{O}_D)$. We assume that $\Sigma$ is a smooth, complete fan and $D=\sum_{i=1}^s a_i D_i$.
We recall that in this case the basis vectors of $H^0(X_\Sigma,\mathcal{O}_D)$ correspond to lattice points inside
$$P_D:=\{m\in M_\q: \langle m,u\rangle \geq -a_u\text{ for all }u\in\Sigma_1\}.$$
In particular, $$\dim H^0(X_\Sigma,\mathcal{O}_D)_m=
\begin{cases}
1\text{ if }m\in P_D\\
0\text{ if }m\not\in P_D.
\end{cases}$$
Referring to Proposition \ref{prop:coh1} we have:
$$m\in P_D\Leftrightarrow \langle m,u\rangle \geq -a_u\text{ for all }u\in\Sigma_1\Leftrightarrow \langle m,u\rangle < -a_u\text{ for none }u\in\Sigma_1\Leftrightarrow$$
$$C_{D,m}\text{ is empty}\Leftrightarrow \tilde H^{-1}(C_{D,m})=1\Leftrightarrow \tilde H^{-1}(C_{D,m})\neq 0.$$
Referring to Proposition \ref{prop:coh1} we have:
$$m\in P_D \Leftrightarrow \langle m,u\rangle +a_u\geq 0\text{ for all }u\in\Sigma_1\Leftrightarrow \Supp((\langle m,u\rangle +a_u)_u)=C_{\{1,\dots,m\}}\Leftrightarrow$$
$$
H^n(\Supp((\langle m,u\rangle +a_u)_u))\neq 0 \Leftrightarrow D+\divi(m)\text{ contributes }1\text{ to the sum in Propostion }
\mbox{\ref{prop:coh2}}
$$
\ex 
Let $X=\p^n$, i.e.~the fan $\Sigma$ has rays $e_1,\dots,e_n,-e_1-\dots-e_n$. 
Let us compute $H^p(X,\mathcal{O}(kD_{e_1})$. (As all divisors $D_{e_i}$ are equivalent this covers all the cases).
We need to consider all decompositions:
$\sum a_i D_i\simeq kD_{e_1}$, which is equivalent to $\sum a_i=k$. Note that if there exists $a_{j_1}\geq 0$ and $a_{j_2}<0$, then $\Supp(a)$ is a nonempty simplex. Hence, all its reduced homology vanish.  

Suppose $k\geq 0$. Then we cannot have all $a_i$ negative, thus we can assume they are all greater or equal to $0$. For each such $a$, the simplicial complex $\Supp(a)$ is an $n$-sphere. Hence, the only nonvanishing homology is the $n$-th one (equal to one), which corresponds to the $0$-th cohomology of the divisor. Thus $H^p(\p^n, \mathcal{O}(kD_{e_1})=0$ for $k\geq 0$ and $p\neq 0$. For $p=0$ we need to count the possible decompositions $\sum_{i=1}^{n+1} a_i=k$ for $a_i\geq 0$. Clearly there are ${k+n}\choose n$ of those. It is worth noting that they naturally correspond to monomials of degree $k$ in $n+1$ variables - the usual description of the basis of sections.

Suppose $k<0$. Now, we only have to consider decompositions with $a_i<0$. In particular, if $-n-1<k$ then all cohomology vanish. As the empty simplicial complex has only $-1$-st reduced homology nonzero (its dimension is equal to $1$ in that case) we see that $H^p(\p^n, \mathcal{O}(kD_{e_1})=0$ for $k<0$ and $p\neq n$. For $p=n$ we need to count the number of decompositions $\sum_{i=1}^{n+1} a_i=k$ with $a_i<0$. Clearly there are $-k-1\choose n$ such possibilities. 
\kex
\ex 
Consider the Hirzebruch surface $X$ with the fan given by $e_1,e_2,-e_1+2e_2,-e_2$. As an example let us compute $H^2(X,\mathcal{O}(-3D_{e_1}-5D_{e_2}))$. As before, (and always) nonvanishing top cohomology corresponds to nonvanishing $-1$-st reduced homology, i.e.~empty simplex. In other words the dimension of the cohomology equals the number of solutions to the diophantine system of equations:
\begin{eqnarray*}
c_1+c_3+2c_4=-5\\
c_2+c_4=-3
\end{eqnarray*}
with $c_i<0$. (This system corresponds to $c_1D_{e_1}+c_2D_{e_2}+c_3D_{-e_1+2e_2}+c_4D_{-e_4}\simeq -5D_{e_1}-3D_{e_2}$.) There are two such solutions.  
\kex
For more relations among simplicial complexes, techniques of counting homology and cohomology of line bundles we refer to \cite[Section 3.2]{lason2011full}.
\section{Gr\"obner basis and triangulations of polytopes}\label{sec:GB}
\dfi[Term order]
An order $<$ on monomials or equivalently on $\N^n$ is called a term order if $0$ is the unique minimal element and $a<b$ implies $a+c<b+c$ for any $a,b,c\in\N^n$. 
\kdfi
\dfi[Initial ideal, Gr\"obner basis]
For a fixed term order $<$ we define $in_<(f)$ to be the unique, largest monomial appearing in $f$ (with nonzero coefficient). For an ideal $I$ we define:
$$in_<(I):=<in_<(f):f\in I>.$$
Caution: it is not enough to take initial forms of (any) generators of $I$ to obtain all generators of $in_<(I)$. 

A set of generators $G$ of $I$ is called \emph{a Gr\"obner basis} (with respect to $<$) if and only if
$$in_<(I):=<in_<(f):f\in G>.$$

A Gr\"obner basis is called minimal if $G$ is \emph{minimal} with respect to inclusion (among sets satisfying equality above). It is called \emph{reduced} if for any $g'\in G$ no monomial in $g'$ is divisible by $in_<(g)$ for some $g\in G$. 

For any ideal $I$ and order $<$ finite Gr\"obner basis exists. Further minimal, reduced Gr\"obner basis is unique (up to scaling).
\kdfi
Most important examples of term orders are \emph{weight orders} $<_\omega$ associated to $\omega\in \r_{\geq 0}^n$. Namely we say that $a\leq b$ if and only if $a\cdot w\leq b\cdot w$. 

(to make it a total order either we need $\omega$ to be irrational or we fix another order $<$ as a tie breaker)

Caution: not every term order is a weight order.

\ob[Proposition 1.11 \cite{Stks}]
For any term order $<$ and any ideal $I$ there exists a non-negative integer vector $\omega \in \N^n$ such that $in_{<\omega}(I)=in_<(I)$.
\kob
\uwa
Geometrically the variety defined by $in_<(I)$ is a degeneration of the variety defined by $I$. In particular, the most important algebraic invariants (like degree or dimension) remain the same. 
This fact is used in algebraic software in order to compute such invariants -- notice that for generic $\omega$ the initial ideal will be monomial.
\kuwa
For $I\subset \c[x_1,\dots,x_n]$ and a fixed term order $<$ we consider a simplicial complex $\Delta_<(I)$ on $\{1,\dots,n\}$ such that $f=\{i_1,\dots,i_k\}$ is a face of $\Delta_<(I)$ if and only if $\prod_{j=1}^k x_{i_j}$ is not in the radical of $in_<(I)$. 
Recall that toric ideals are in variables that have natural interpretation as lattice points. Suppose that our toric variety is defined by a lattice polytope $P$. Our plan is to realize $\Delta_<(I)$ as a 'subdivision' of $P$. Let $P$ be a $d-1$ dimensional lattice polytope in $\r^{d-1}\times \{1\}\subset\r^d$.
\dfi[Regular triangulation]\label{def:regtrian}
Given a subset of lattice points $P=\{p_1,\dots,p_n\}$ and $\omega \in \r^n$ we define a subdivision $\Delta_\omega$ by:

$f=\{p_{i_1},\dots,p_{i_k}\}$ is a face of $\Delta_\omega$ if any only if there exists $c\in \r^d$ such that:
$$p_j\cdot c=\omega_j\text{ if }j\in\{i_1,\dots,i_k\}\text{ and }$$
$$p_j\cdot c<\omega_j\text{ if }j\not\in \{i_1,\dots,i_k\}.$$
A generic $\omega$ defines a triangulation $\Delta_\omega$ and any such triangulation is called \emph{regular}.
\kdfi
\uwa
One should imagine a two dimensional polytope $P$ and the heights $w_i$ above $p_i$ as points in the third dimension. The linear forms $c$ are obtained by putting a sheet of paper from below and stopping at some heights. (If $\omega$ is generic such a sheet of paper will touch exactly three points. These are the triangles (after projecting back to $P$) of the triangulation $\Delta_\omega$.) 
\kuwa
\uwa
Not any triangulation is regular \cite[Example 8.2]{Stks}.
\kuwa
\dfi[Secondary fan]
The regular triangulations in Definition \ref{def:regtrian} depend on $\omega$. However, the set of $\omega \in \r^n$ that give rise to the same regular triangulation is a cone. These cone form a complete fan in $\r^n$ called the \emph{secondary fan}.
\kdfi
\tw[Theorem 8.3 \cite{Stks}]
$\Delta_{<_\omega}(I_P)=\Delta_\omega$
\ktw
\wn
$$Rad(in_\omega(I_P))=<x^b:b\text{ is a minimal nonface of }\Delta_\omega>=\bigcap_{\sigma\in\Delta_\omega}<x_i:i\not\in\sigma>$$
\kwn
What is the multiplicity of a prime ideal in $in_\omega(I_P)$?
\tw[Theorem 8.8 \cite{Stks}]
For $\sigma \in\Delta_\omega$ of dimension $d-1$ the normalized volume of $\sigma$ equals the multiplicity of the prime ideal $<x_i:i\not\in\sigma>$ in $in_\omega(I_P)$. 
\ktw
\wn[Corollary 8.9 \cite{Stks}]
The initial ideal $in_\omega(I_P)$ is square-free if and only if the corresponding regular triangulation is unimodular.
\kwn

\section{Cuts and Splits}\label{sec:cutsandsplits}
Let $G=(V,E)$ be a graph. Following \cite{Cutsandsplits} we consider two polynomial rings:
$$\C[q]:=\C[q_{A|B}: A\cup B=V, A\cap B=\emptyset],$$
$$\C[s,t]:=\C[s_{ij},t_{ij}: \{i,j\}\in E].$$
In the first ring the variables correspond to partitions of $V$, in the second to each edge we associate two variables. For a given partition $A|B$ we define a subset $\Cut(A|B)\subset E$ of cut edges by:
$$\Cut(A|B):=\{\{i,j\}\in E: i\in A, j\in B\text{ or }j\in A, i\in B\}.$$
There is a natural monomial map:
$$\C[q]\rightarrow \C[s,t], \quad q_{A|B}\rightarrow \prod_{\{i,j\}\in \Cut(A|B)} s_{ij}\prod_{\{i,j\}\in E\setminus\Cut(A|B)}t_{ij}.$$
As the map is monomial its kernel $I_G$ is a prime binomial ideal, and as the monomials are homogeneous, so is $I_G$. 

Sturmfels and Sullivant posed several conjectures how combinatorics of $G$ relates to algebraic properties of $\C[q]/I_G$. In spite of progress on this topic \cite{Engstrom} many still remain open.
\hip[\cite{Cutsandsplits}, Conjecture 3.7]\label{hip:cutsNormal}
$\C[q]/I_G$ is normal if and only if $G$ is free of $K_5$ minors.
\khip
By \cite{OhsugiCuts} we know that the class of graphs for which $\C[q]/I_G$ is normal is minor closed. Moreover, $\C[q]/I_{K_5}$ is not normal. We will now show how Conjecture \ref{hip:cutsNormal} implies the famous four color theorem. Originally the idea is due to  David Speyer. First we need some results on the polytope representing the toric variety. It is called the cut polytope. Its vertices can be described as indicator vectors on the set of edges of $G$, corresponding to cut edges, for any partition of vertices. From now one we assume that $G$ is a graph without $K_5$ minor.
We have the following theorem due to Seymour.
\tw[\cite{SeymourMatroids}, \cite{CutPolytope}, Corollary 3.10]\label{tw:PGfacets}
Let $P_G\subset \{1\}\times \q^{|E|}\subset \q^{|E|+1}$ be the cut polytope. For any edge $e\in E$ we denote by $x_e$ the corresponding coordinate in $\q^{|E|}$ and by $x_0$ the first coordinate.
The cone over the polytope $P_G$ is defined by the following inequalities:
$$0\leq x_e\leq x_0,\text{ for each edge }e\text{ that does not belong to a triangle},$$
$$\sum_{e\in F}x_e\leq (|F|-1)x_0+\sum_{e\in C\setminus F}x_e,\text{ for any cordless cycle }C, F\subset C, |F|\text{ odd}.$$
\ktw
\lem\label{lem:point23}
The point $p=(3,2,\dots,2)\in \Z^{|E|+1}$ belongs to the lattice spanned by $P_G$ and to the cone over this polytope.
\klem
\dow 
Summing up all points corresponding to partitions $\{v\}|V\setminus\{v\}$ we obtain the point $(|V|,2,\dots,2)$. As the partition $\emptyset|V$ corresponds to $(1,0,\dots,0)$ we see that $p$ indeed is in the lattice spanned by $P_G$. It is straightforward to check that $p$ satisfies the inequalities in Theorem \ref{tw:PGfacets}.
\kdow
\ob\label{prop:con=>4col}
Conjecture \ref{hip:cutsNormal} implies the four color theorem.
\kob
\dow 
Consider $p$ from Lemma \ref{lem:point23}. If $P_G$ is normal, then $p=p_1+p_2+p_3$, where $p_i$ corresponds to a partition $A_i|B_i$. Hence, we have three partitions, such that any edge belongs to precisely two of them. We define four subsets of $V$:
\begin{enumerate}
\item $(A_1\cap A_2\cap A_3)\cup(B_1\cap B_2\cap B_3)$,
\item $(A_1\cap A_2\cap B_3)\cup(B_1\cap B_2\cap A_3)$,
\item $(A_1\cap B_2\cap A_3)\cup(B_1\cap A_2\cap B_3)$,
\item $(B_1\cap A_2\cap A_3)\cup(A_1\cap B_2\cap B_3)$.
\end{enumerate}
Clearly, these subsets define a partition of $V$, i.e.~are disjoint and every vertex belongs to one of them. We now prove that this is a proper coloring. As the choice between $A_i$ and $B_i$ was arbitrary, it is enough to prove there are no edges among vertices in $(A_1\cap A_2\cap A_3)\cup(B_1\cap B_2\cap B_3)$. Indeed, if both vertices of such edge $e$ belonged to $A_1\cap A_2\cap A_3$ or $B_1\cap B_2\cap B_3$ then $x_e(p)=0$. However, if one vertex belongs to $A_1\cap A_2\cap A_3$ and the other to $B_1\cap B_2\cap B_3$ then $x_e(p)=3$. This finishes the proof.
\kdow
\uwa
We showed that Conjecture \ref{hip:cutsNormal} implies four colorability of any graph without $K_5$ minor. This apparently stronger statement was in fact classically known to be equivalent to the four color theorem and is a special case of a more general Hadwiger conjecture.
\kuwa
\uwa 
In the same way one can show that a $4$-coloring of $G$ induces a decomposition of $p=p_1+p_2+p_3$. The three partitions come from dividing the four colors into two groups of two colors. 
\kuwa

\section{Toric varieties and matroids}\label{sec:TVandMat}

Let $M$ be a matroid on a ground set $E$ with the set of bases $\mathfrak{B}\subset\mathcal{P}(E)$ (the reader is referred to \cite{Ox} for background of matroid theory). Let
$S_M:=\C[y_B:B\in\mathfrak{B}]$
be a polynomial ring. Let $\varphi_M$ be the $\C$-homomorphism:
$$\varphi_M:S_M\ni y_B\rightarrow\prod_{e\in B}x_e\in\C[x_e:e\in E].$$
The \emph{toric ideal of a matroid} $M$, denoted by $I_M$, is the kernel of the map $\varphi_M$. 
\tw[\cite{white1977basis}]\label{thm:matroidnormal}
The ideal $I_M$ defines a \emph{normal} toric variety.
\ktw
We recall the following theorem.
\tw[\cite{Stks} Theorem 13.14]\label{tw:deggen}
The toric ideal $I_P$ associated to a normal polytope $P$ is generated in degree at most $\dim P$.  
\ktw
As a corollary of Theorem \ref{thm:matroidnormal} and \ref{tw:deggen} we obtain the following.
\wn\label{cor:degE}
For any matroid $M$ on the ground set $E$ the toric ideal $I_E$ is generated in degree at most $|E|$.
\kwn

A result of Gijswijt and Regts strengthens Theorem \ref{thm:matroidnormal}.

\tw[\cite{GiRe12}]\label{thm:matroidnormal2}
(Poly)Matroid base polytope (associated to the monomials defining $\phi_M$) has an Integer Caratheodory Property. That is, if $P_M$ is a matroid base polytope, then every integer vector in $kP_M$ is a positive, integral sum of affinely independent integer vectors from $P_M$ with coefficients summing up to $k$. In particular, Caratheodory rank of $P_M$ is as low as possible, it is equal to the dimension of $P_M$ plus $1$. 
\ktw

\subsection{Representable matroid}

This subsection is not needed in what follows. However, it provides additional motivation to study $V(I_M)$ from the point of view of algebraic geometry.
Let $M$ be a representable matroid realized by vectors $v_1,\dots,v_k$ spanning a $d$-dimensional vector space $V$. 
We have a natural map $\C^k\ni e_i\rightarrow v_i\in V$ with kernel $K\in G(k-d,k)$. Note, that on $\C^k$ acts a $k$-dimensional torus $T$ inducing an action on the Grassmannian $G(k-d,k)$.
\tw\cite{GGMS}\label{thm:orbitsinG}
The toric variety $V(I_{M^*})$ is isomorphic to the closure of the orbit $G\cdot K$ in $G(k-d,k)$.
\ktw
\dow 
Let $f_1,\dots,f_{k-d}$ be the basis of $K$.  Let $M_k$ be a $(k-d)\times k$  matrix with $i$-th row corresponding to $f_i$. A given Plu\"ucker coordinate $(e_{a_1}\wedge\dots\wedge e_{a_{k-d}})^*(f_1\wedge\dots\wedge f_{k-d})$ of $K$ equals the maximal minor of $M_k$ distinguished by the columns $\{a_1,\dots,a_{k-d}\}$. Further, the torus $T$ acts on $e_{a_1}\wedge\dots\wedge e_{a_{k-d}}$ with weight corresponding to a lattice point in $\Z^k$ that has coordinates indexed by $a_i$ equal to $1$ and all other equal to $0$. Hence: 
\begin{itemize}
\item the polytope $P_{M^*}$ has lattice points corresponding to indicator vectors of complements of bases of $M$,
\item the polytope of the toric variety $\overline{T\cdot K}$ has lattice points corresponding to indicator vectors of subsets of $k-d$ columns of $M_k$ giving a nonzero minor.
\end{itemize}
Thus, it remains to show that a given minor of $M_k$ is nonzero if an only if the corresponding vectors $v_i$ form a complement of a basis of $V$. Fix a set $S=\{a_1,\dots,a_{k-d}\}\subset \{1,\dots,k\}$ and denote its complement by $S'$.
We have following equivalences:

$\{v_i\}_{i\in S}$ is a complement of a basis of $V$ $\Leftrightarrow$ $\{v_i\}_{i\in S'}$ is a basis of $V$ $\Leftrightarrow$ $\{f_j\}_{j=1,\dots,k-d}\cup \{e_i\}_{i\in S'}$ is a basis of $\C^k$ $\Leftrightarrow$ the minor of $M_k$ distinguished by $S$ is nonzero.
\kdow
\uwa
If we fix a basis of the space $V$ we obtain a matrix representation of $v_i$ which distinguishes a point in $G(d,k)$. The closure of the torus orbit of this point is isomorphic to $V(I_M)$.
\kuwa

\subsection{Open problems}\label{sec:matOpen}

Let us present the major open problems concerning toric ideals associated to matroids.
\hip[Weak White's conjecture \cite{White}]\label{hip:Whiteweak}
For any matroid $M$ the ideal $I_M$ is generated by quadrics.
\khip
We say that $y_{B_1}y_{B_2}-y_{B_3}y_{B_4}$ for $B_i\in \mathfrak{B}$ is a \emph{symmetric exchange} if $B_3=(B_1\setminus\{b_1\})\cup\{b_2\}$ and $B_4=(B_2\setminus\{b_2\})\cup\{b_1\}$ for some $b_1\in B_1\setminus B_2$, $b_2\in B_2\setminus B_1$ (see \cite{La15} for other exchange properties).
\hip[White's conjecture \cite{White}]\label{hip:White}
For any matroid $M$ the ideal $I_M$ is generated by symmetric exchanges.
\khip
However, even the following is open.
\hip\label{hip:quadricsbysym}
Every quadric in $I_M$ is a linear combination of symmetric exchanges.
\khip
In view of Theorem \ref{thm:orbitsinG} the following turns out to be a very important open problem, that is weaker than White's conjecture.
\hip\label{hip:realizable}
For any representable matroid $M$, the ideal $I_M$ is generated by quadrics.
\khip

Herzog and Hibi write, that they `cannot escape from the temptation' to ask the following, stronger questions.

\begin{question}[\cite{HerzogHibi}]\label{qu:HerzogHibi}
Let $M$ be a discrete (poly)matroid.
\begin{enumerate}
\item Does the toric ideal $I_M$ possess a quadratic Gr\"{o}bner basis?
\item Is the base ring $\C[y_B:B\in\mathfrak{B}]/I_M$ Koszul?
\end{enumerate}
\end{question}

A positive answer to the first question would imply the following conjecture, that is a strengthening of Theorem \ref{thm:matroidnormal2}.
\hip
For any matroid $M$ the basis polytope $P_M$  has a unimodular triangulation.
\khip

A \emph{$k$-matroid} is a matroid whose ground set can be partitioned into $k$ pairwise disjoint bases. We call a basis of a $k$-matroid \emph{complementary} if its complement can be partitioned into $k-1$ pairwise disjoint bases. The \emph{basis graph} of a matroid is a graph with vertices corresponding to bases and an edge between two bases that differ by a pair of elements. 
The \emph{complementary basis graph} of a $k$-matroid is the restriction of its basis graph to complementary bases.

\hip\cite{La16}\label{ConjectureConnected}
Complementary basis graph of a $k$-matroid is connected.
\khip

Notice that Conjecture \ref{ConjectureConnected} for $k=2$ coincides with Conjecture \ref{hip:quadricsbysym} in a non-commutative setting.
\uwa
For a $2$-matroid it is not known even if some antipodal bases $B_1,B_2$ (that is, bases such that $B_1\sqcup B_2=E$) are always connected in the complementary basis graph. Positive answer would imply that commutative and non-commutative settings of Conjecture \ref{hip:quadricsbysym} are equivalent.
\kuwa
\hip\cite{La16}\label{Conjecturekr+1}
Let $k\geq 2$, and let $M$ be a matroid of rank $r$ on the ground set $E$ of size $kr+1$. Suppose $x,y\in E$ are two elements such that both sets $E\setminus x$ and $E\setminus y$ can be partitioned into $k$ pairwise disjoint bases. Then there exist partitions of $E\setminus x$ and $E\setminus y$ into $k$ pairwise disjoint bases which share a common basis.
\khip

If $k\geq 2^{r-1}+1$, then the above Conjecture \ref{Conjecturekr+1} holds \cite{La16}.

\ob\cite{La16}
White's Conjecture \ref{hip:White} implies Conjecture \ref{ConjectureConnected}.\newline
Conjunction of Conjectures \ref{ConjectureConnected} and \ref{Conjecturekr+1} implies White's Conjecture \ref{hip:White}.
\kob

\subsection{Known facts}

Several special cases of Conjecture \ref{hip:White} are known.
\tw 
Conjecture \ref{hip:White} holds for graphic \cite{Blasiak}, sparse paving \cite{Bonin}, strongly base orderable \cite{jaMichal}, and rank $\leq 3$ \cite{Ka10} matroids.

Further, the classes of matroids for which the toric ideal is generated by
quadrics and that has quadratic Gr\"obner bases, is closed under series and
parallel extensions, series and parallel connections, and 2-sums \cite{shibata2016toric}.
\ktw
Further, partial results are known in general. Let $J_M\subset I_M$ be the ideal generated by symmetric exchanges. Let $m\subset S_M$ be the irrelevant ideal, i.e.~the maximal ideal generated by all the variables.
\tw\cite{jaMichal}\label{tw:jaMichalMatroid}
For any matroid $M$, we have $J_M:m^\infty = I_M$. That is, $\Proj S_M/J_M$ is equal to $\Proj S_M/I_M$. In particular, Conjecture \ref{hip:White} holds on set-theoretic level.
\ktw

We can rephrase the above by saying that homogeneous components of the ideals $I_M$ and $J_M$ of a matroid $M$, are equal starting from some degree $f(M)$. Even more is known.

\tw\cite{La16}
Homogeneous components of ideals $I_M,J_M$ of a matroid $M$ of rank $r$, are equal starting from degree $f(r)$ depending only on $r$. 
\ktw

It follows that checking if White's Conjecture \ref{hip:White} is true for matroids of a fixed rank is a decidable (finite) problem.

The first question of Herzog and Hibi \ref{qu:HerzogHibi} is already difficult for special classes of matroids. For uniform matroids, it was proved by Sturmfels \cite{Stks}. Later, it was extended to base-sortable matroids by Blum \cite{Bl01} and further by Herzog, Hibi and Vladoiu \cite{HHV}. Schweig \cite{Sc11} proved it for lattice path matroids. The only general result concerning Gr\"{o}bner bases of toric ideals of matroids is the following.

\tw\cite{La17}
The toric ideal of a matroid of rank $r$ posesses a Gr\"{o}bner basis of degree at most $(r+4)!$.
\ktw

Conca \cite{Conca} proved that the answer to the second Question \ref{qu:HerzogHibi} is positive for transversal polymatroids.

\subsection{Matroids and finite characteristics}\label{sec:MatFC}
The original idea to look at implications of Theorem \ref{tw:jaMichalMatroid} in case of finite characteristics and pureness is due to Matteo Varbaro.
\dfi[Pure, F-pure \cite{HochsterRoberts}]\label{def:pure}
A morphism of rings $R\rightarrow S$ is called \emph{pure} if for any $R$-module $M$ the map $M\ni m\rightarrow m\otimes 1\in M\otimes_R S$ is injective.

If $R$ is an algebra over a field $k$ of characteristic $p$ we say that $R$ is F-pure if the Frobenius morphism $R\ni r\rightarrow r^p\in R$ is pure. 
\kdfi
\lem\label{lem:pure=>reduced}
If a $k$-algebra $R$ is F-pure, then it is reduced.
\klem
\dow 
Take $M=R$ in Definition \ref{def:pure}. For any $a\in \N$ and for any $r\in R$ we have:
$$r^a\rightarrow r^a\otimes 1=1\otimes r^{pa}.$$
Thus, we have $r^a\neq 0\Rightarrow r^{pa}\neq 0$, which proves the lemma.
\kdow
\ob\label{ob:ConeFpure}
Let $C$ be a cone in a lattice $M$. Then $k[C]$ is F-pure.
\kob
\dow 
The proof is a combination of the following three facts:
\begin{enumerate}
\item A (Lauarent) polynomial ring over $k$ is F-pure (more generally a Noetherian regular $k$-algebra) \cite[Proposition 5.14]{HochsterRoberts}.
\item 
If $S$ if F-pure and $S=R\oplus M$ as $R$ modules for some $R$ module $M$ then $R$ is $F$ pure \cite[Proposition 1.3]{Fedder}.
\item One can realize $C=\Z^r\times C'$, where $C'=\Z_+^m\cap H$ and $H$ is a linear subspace. In particular, the ring $\C[C']$ as a module is a factor of $\C[\Z_+^m]$ \cite[p.~63]{BG}.
\end{enumerate}
\kdow
\tw 
For any matroid $M$ the following are equivalent:
\begin{enumerate}
\item $I_M$ is generated by symmetric exchanges (White's conjecture holds),
\item $J_M$ is saturated with respect to $m$,
\item $J_M$ is prime,
\item $J_M$ is primary,
\item $S_M/J_M$ is an integral domain,
\item $S_M/J_M$ is reduced,
\item The localisation $(S_M/J_M)_m$ is reduced,
\item $S_M/J_M$ is F-pure over some field $k$ of finite characteristic,
\item $S_M/J_M$ is F-pure over any field $k$ of finite characteristic,
\item The localisation $(S_M/J_M)_m$ is F-pure over some/any field.
\end{enumerate}
\ktw
\dow 
Equivalences $i)-vii)$ are direct consequences of Theorem \ref{tw:jaMichalMatroid}. 

For a fixed field $k$ point $i)$ implies $viii)$ by a combination of Theorem \ref{tw:normalpoly} and Proposition \ref{ob:ConeFpure} (recall that a toric algebra is normal if the corresponding monoid is a cone - Theorem \ref{tw:charnormal}). 

The equivalence of local and global case follows from \cite[Proposition 1.3]{Fedder}.

Point $viii)$ implies $vi)$ by Lemma \ref{lem:pure=>reduced}.

It remains to show that $viii)$ is equivalent to $ix)$ or in other words that White's conjecture does not depend on the underlying field. More general statement is explained in Lemma \ref{lem:fieldindependent}. 
\kdow
\lem\label{lem:fieldindependent}
Consider a finite set of characters $S\subset M\simeq \Z^n$ generating a monoid $\tilde S$. Suppose for some field $k$ the kernel $I_k$ of the map:
$$k[x_1,\dots,x_n]\rightarrow k[\tilde S]$$
is generated by binomials $B=\{m_1-m_1',\dots,m_l-m_l'\}$. Then the same binomials generate the toric ideal $I_{k'}$ associated to $S$ over any other field $k'$.
\klem
\dow 
Reasoning as in Theorem \ref{tw:binomials} we see that the binomials in the toric ideal do not depend on the field and simply correspond to integral relations among lattice points in $S$. Further, the ideal is always spanned by such binomials. Thus we only need to show the following:

Any binomial $m-m'\in I_{k'}$ is equal to $\sum_j n_j(m_{i_j}-m_{i_j}')$, where $n_j$ are monomials and $m_{i_j}-m_{i_j}'\in B$.

We know that $m-m'=\sum \lambda_j m_j (m_{i_j}-m_{i_j}')$ for some $\lambda_j\in k^*$. By definition, say that the monomial $m_jm_{i_j}$ is equivalent to $m_jm_{i_j}'$. We generate the equivalence relation and we want to prove that $m$ is equivalent to $m'$. We may subdivide the terms in $\sum \lambda_j m_j (m_{i_j}-m_{i_j}')$ according to the equivalence class they belong to. Consider all terms in the same class as $m$. This is a sum of monomials with some coefficients, one of which is $m$ with coefficient one. However, the sum of these coefficients needs to be $0$ and the only other monomial with nonzero coefficient that may appear is $m'$. 
\kdow
Let us restate a criterion of F-pureness due to Fedder. For an ideal $I$ we let $I^{[p]}$ be the ideal generated by $\{i^p:i\in I\}$.
\tw[Proposition 1.3 and 1.7 \cite{Fedder}]
Let $J$ be a homogeneous ideal in a polynomial ring $S$ over $k$ of characteristic $p$. Then $S/J$ is F-pure if and only if $J^{[p]}:J\not\subset m^{[p]}$.
\ktw 
\wn
White's conjecture \ref{hip:White} holds if and only if there exists $f\not \in m^{[p]}$ such that for any symmetric exchange $g$ we have $gf\in J_M^{[p]}$. Further, we may (but do not have to) assume that we work over the field $\Z_2$, i.e.~$p=2$ and $f$ is homogeneous of degree equal to the number of bases of the matroid. In such a case $f\not \in m^{[p]}$ translates to $f$ being a sum of monomials one of which is a (squarefree) product of all variables in $S_M$.
\kwn
We note that in case $k=\Z_2$ the polynomial $f$ distinguishes some multisets of bases - it would be great to understand combinatorial meaning of those. In each particular case of $M$ we can compute (examples of) $f$.
\ex\label{exm:Fedder}
The grading below encodes a uniform rank two matroid on a ground set with four elements. First we compute the toric ideal in Macaulay2.
\begin{verbatim}
loadPackage,"Normaliz"
loadPackage("MonomialAlgebras",Configuration=>{"Use4ti2"=>true})
L={{1,1,0,0},{1,0,1,0},{1,0,0,1},{0,1,1,0},{0,1,0,1},{0,0,1,1}};
d=#L;
R=ZZ/2[a_1..a_d,Degrees=>L]
J=binomialIdeal R; m=ideal(gens R);
\end{verbatim}
We now check Fedder's criterion.
\begin{verbatim}
Jp=0; mp=0;
for i from 0 to ((numgens J)-1) do {Jp=Jp+ideal( ((gens J)_i_0)^2)}
for i from 0 to ((numgens m)-1) do {mp=mp+ideal( ((gens m)_i_0)^2)}
isSubset(Jp:J,mp)
Q=R/mp;
Red=sub(Jp:J,Q)
\end{verbatim}
The answer the program gives is:
$$f={a}_{2} {a}_{3} {a}_{4} {a}_{5}+{a}_{1} {a}_{3} {a}_{4} {a}_{6}+{a}_{1} {a}_{2} {a}_{5} {a}_{6}.$$

Below we present the example of a graphic matroid corresponding to a square with one diagonal.
\begin{verbatim}
L={{1,1,1,0,0},{1,1,0,1,0},{1,0,1,1,0},{0,1,1,1,0},{1,0,1,0,1},{0,1,0,1,1},
{1,0,0,1,1},{0,1,1,0,1}};
d=#L;
R=ZZ/2[a_1..a_d,Degrees=>L]
J=binomialIdeal R; m=ideal(gens R);
Jp=0; mp=0;
for i from 0 to ((numgens J)-1) do {Jp=Jp+ideal( ((gens J)_i_0)^2)}
for i from 0 to ((numgens m)-1) do {mp=mp+ideal( ((gens m)_i_0)^2)}
isSubset(Jp:J,mp)
Q=R/mp;
Red=sub(Jp:J,Q)
\end{verbatim}
Here:
$$f={a}_{1} {a}_{4} {a}_{5} {a}_{6} {a}_{7}+{a}_{2} {a}_{3} {a}_{5} {a}_{6} {a}_{8}+{a}_{2} {a}_{4} {a}_{5} {a}_{7} {a}_{8}+{a}_{
1}
      {a}_{3} {a}_{6} {a}_{7} {a}_{8}.$$
\kex
One of the problems that we encountered is that $f$ may be not uniquely specified (e.g. in Example \ref{exm:Fedder} one can multiply $f$ by a variable). Further, it is hard to see a general pattern for $f$. 
\section{Toric varieties and phylogenetics}\label{sec:TVandPhylo}
In this section we present a construction of a family of toric varieties. It is inspired by phylogenetics - a science that aims at reconstruction of the history of evolution. The basic statistical picture is as follows. We start from a (usually large) family of parameters that correspond to various probabilities of mutations (and probability distribution of the common ancestor). These parameters are unknown. However, if we knew them then we could answer questions of type:

'(On a given position in DNA string) what is the probability that a human has $C$, a gorilla has $C$ and a guenon has $A$?'

Here, we represent DNA as strings of characters: $A,C,G,T$. Thus, we obtain a map $m$ from the parameter space to the space of joint probability distributions of states of species that we consider. The latter is a huge space! The states are indexed by a choice of a letter for each species we consider - the dimension is $4$ to the power equal to the number of species. Note that there is a distinguished point $P$ in this space: going through whole DNA sequences (that are very long) biologists and statisticians can count the number of times they encounter $C$ for human, $C$ for gorilla and $A$ for guenon etc. In other words, we are getting a(n approximation) of the probability distribution on the (joint) states of the species. 

What are the main questions that we would like to answer?

\begin{itemize}
\item Does the point $P$ belong to the image of $m$? (If the answer is no, then it means that some of our assumptions are wrong in the statistical model we have chosen.)
\item Can we identify the parameters $m^{-1}(P)$?
\end{itemize}
Let us focus on the first question. First, we notice that in most interesting examples the map $m$ is algebraic.
The approach algebraic geometry proposes is to consider the defining equations of the Zariski closure $X$ of the image of $m$. 

The algebraic variety $X$ we obtain depends on the statistical model we choose, i.e.~what we assume on the parameters. The most 'universal' is the General Markov Model that leads to secant varieties of Segre products. Let us discus a different class of so-called group-based models. Their biological motivation comes from the observation that there are certain symmetries among probabilities of mutations. These symmetries can be encoded by the group action. For example the famous Kimura $3$-parameter model \cite{3Kimura} relies on the fact that $\{A,C,G,T\}$ is naturally divided into two subsets: purines $\{A,G\}$ and pyrimidines $\{C,T\}$. In mathematical language the group $\Z_2\times \Z_2$ naturally acts on $\{A,C,G,T\}$ and this influences the geometry of the variety $X$ for the $3$-Kimura model.
\subsection{Construction for group-based models} In general to determine the variety $X$ we need a tree $T$ (that determines how species mutated) and a model (in our case determined by a finite abelian group $G$). Here, we do not present the general construction referring to \cite{SS, 4aut, jaDissert}. On the other hand, we describe in detail the toric structure of the variety $X$ in case when it is a so-called star or claw-tree $K_{1,n}$, i.e.~a tree with one inner vertex and $n$ leaves. The general case may be also obtained by \emph{toric fiber product} \cite{Sethtfp}. As we will see many conjectures address the case of $K_{1,n}$.

\begin{df}[{\bf Flow} \cite{JaAdvGeom}, \cite{BW}]  \label{def:flow}
Let $G$ be a finite abelian group and $n\in \mathbb N$. A flow 
is a sequence of $n$ elements of $G$ summing up to $0\in G$, the neutral element of $G$. The set of flows is equipped with a group structure via the coordinatewise action. 
The group of flows $\Gf$ is (non-canonically) isomorphic to $G^{n-1}$.
\end{df}
\begin{df}[{\bf Polytope $P_{G,n}$}, \cite{jaEmanuele}, \cite{SS}]\label{def:PGn}
Consider the lattice $M\cong \Z^{|G|}$ with a basis corresponding to elements of $G$. Consider $M^n$ with the basis $e_{(i,g)}$ indexed by pairs $(i,g)\in [n]\times G$. We define an injective map of sets:
$\Gf\rightarrow M^n,$
by $(g_1,\dots,g_n)\longmapsto \sum_{i=1}^n e_{(i,g_i)}$. The image of this map defines the vertices of the polytope $P_{G,n}$.
\end{df}

\begin{exm}[{\bf \cite{jaZ3nowy}}]

For $G=(\Z_2,+)$ and $n=3$, we have four flows:
$$(0,0,0),(0,1,1),(1,0,1),(1,1,0)\in \Z_2\times \Z_2\times \Z_2.$$

\noindent Hence, the polytope $P_{\Z_2,3}$ has the following four vertices corresponding to the flows above:
$$(1,0,1,0,1,0),(1,0,0,1,0,1),(0,1,1,0,0,1),(0,1,0,1,1,0)\in\Z^2\times \Z^2\times \Z^2,$$
where $(1,0)\in \Z^2$ corresponds to $0\in \Z_2$ and $(0,1)\in\Z^2$ corresponds to $1\in \Z_2$.
\end{exm}

A more sophisticated example is presented in \cite[Example 4.1]{JaJalg}. It turns out that the phylogenetic variety $X$ - for group based models - is toric and corresponds to the polytope $P_{G,n}$.
We already know by Theorem \ref{tw:binomials} that binomials in the toric ideal correspond to integral relations among lattice points. However, for group based models it is easier to work with flows. 
Binomials may be identified with a pair of tables of the same size $T_0$ and $T_1$ of elements of $G$, regarded up to row permutation. Each row of such tables has to be a flow. The identification is as follows. Every binomial is a pair of monomials; the variables in such monomials correspond to flows, given by a collection of $n$ elements in $G$. Every monomial is viewed as a table, whose rows are the variables appearing in the monomial; the number of rows of the corresponding table is the degree of the monomial. Consequently, a binomial is identified with the pair of tables encoding the two monomials respectively. \\
\indent For a finite abelian group $G$ and the graph $K_{1,n}$ the associated toric variety (represented by the polytope $P_{G,n}$) will be denoted by $X(G,K_{1,n})$. A binomial belongs to $I(X(G , K_{1,n}))$ if and only if the two tables are \emph{compatible}, i.e.~for each $i$, the $i$-th column of $T_0$ and the $i$-th column of $T_1$ are equal as multisets.\\ 
\indent In order to generate a binomial -- represented by a pair of tables $T_0$, $T_1$ -- by binomials of degree at most $d$ we are allowed to select a subset of rows in $T_0$ of cardinality at most $d$ and replace it with a compatible set of rows, repeating this procedure until both tables are equal. 

\begin{exm}[{\bf \cite{jaZ3nowy}}]\label{ex:rel}

For $G=(\Z_2,+)$ and $n=6$ consider the following two compatible tables:
$$
T_0=\begin{bmatrix}
\color{red}{1} & \color{red}{1} & \color{red}{1} &  \color{red}{1} & \color{red}{1} & \color{red}{1} \\
\color{red}{0} & \color{red}{0} & \color{red}{0} & \color{red}{0} & \color{red}{0} & \color{red}{0} \\
\color{brown}{1} & \color{brown}{1} & \color{brown}{0} & \color{brown}{0} & \color{brown}{0} & \color{brown}{0}\\
\end{bmatrix}
\textnormal{ and   } 
T_1=\begin{bmatrix}
\color{blue}{0} & \color{blue}{1} & \color{blue}{0} & \color{blue}{1} & \color{blue}{0}& \color{blue}{0} \\
1 & 0 & 1 & 0 & 0 & 0 \\ 
1 & 1 & 0 & 0 & 1 & 1 \\
\end{bmatrix}.
$$
\noindent 
Note that the red subtable of $T_0$ is compatible with the table 
$$
T'=\begin{bmatrix}
\color{blue}{0} & \color{blue}{1} & \color{blue}{0} & \color{blue}{1} & \color{blue}{0}& \color{blue}{0} \\
\color{brown}{1} & \color{brown}{0} & \color{brown}{1} & \color{brown}{0} & \color{brown}{1} & \color{brown}{1}\\
\end{bmatrix}.
$$
\noindent Hence, we may \emph{exchange} them obtaining:
$$
\tilde T_0=\begin{bmatrix}
\color{blue}{0} & \color{blue}{1} & \color{blue}{0} & \color{blue}{1} & \color{blue}{0}& \color{blue}{0} \\
\color{brown}{1} & \color{brown}{0} & \color{brown}{1} & \color{brown}{0} & \color{brown}{1} & \color{brown}{1}\\ 
\color{brown}{1} & \color{brown}{1} & \color{brown}{0} & \color{brown}{0} & \color{brown}{0} & \color{brown}{0}\\
\end{bmatrix}.
$$

\noindent Note that $T_0$ and $\tilde T_0$ are \emph{compatible}.  Now, the brown subtable of $\tilde T_0$ is compatible with the table
$$
T''=\begin{bmatrix}
1 & 0 & 1 & 0 & 0 & 0 \\
1 & 1 & 0 & 0 & 1 & 1 \\
\end{bmatrix}.
$$
\noindent Finally, we exchange them obtaining $T_1$. Hence we have a sequence of tables $T_0 \rightsquigarrow \tilde T_0 \rightsquigarrow T_1$. More specifically, we started from a degree three binomial given by the pair $T_0, T_1$ and we generated it using degree two binomials of degree two. 

\end{exm}
\begin{df}[{\bf Phylogenetic complexity} \cite{SS}]

Let $K_{1,n}$ be the star with $n$ leaves, and let $\phi(G,n)$ be the maximal degree of a generator in a minimal generating set of $I(X(G,K_{1,n}))$. We define the phylogenetic complexity $\phi(G)$ of $G$ to be $\sup_{n\in \mathbb N} \phi(G,n)$. 

\end{df}

A new package to deal with phylogenetic group-based models appeared recently \cite{banos2016phylogenetic} for Macaulay2. The software to generate polytopes $P_{G,n}$ is presented in \cite{DBM}.
\subsection{Further properties of group-based models}\label{sec:GBM}
In studying group-based models, Buczy{\'n}ska and Wie{\'s}niewksi \cite{BW}, \cite{MR2892983} made the startling observation that in the case $G = \Z/2\Z$ the Hilbert function of the affine semigroup algebra $\C[M_{\Gamma,G}]$ associated to a graph $\Gamma$ (with respect to an appropriate grading) only depends on the number of leaves $n$ and the first Betti number $g$ of $\Gamma$. The explanation for this phenomenon was provided by Sturmfels and Xu \cite{SX} and Manon \cite{M}, where it was shown that the phylogenetic statistical models $M_{\Gamma, \Z/2\Z}$ are closely related to the Wess-Zumino-Witten (WZW) model of conformal field theory, and the moduli space $\mathcal{M}_{C, \vec{p}}(SL_2(\kk))$ of rank $2$ vector bundles on an $n$-marked algebraic curve $(C, \vec{p})$ of genus $g$. In particular, the total coordinate ring $\mathcal{V}_{C, \vec{p}}(SL_2(\kk))$ of this space, which is known to be a direct sum of the so-called conformal blocks of the WZW model, is shown to carry a flat degeneration to each affine semigroup algebra $\kk[M_{\Gamma, \Z/2\Z}]$. Flat degeneration preserves Hilbert polynomials, explaining the coincidence among the $\kk[M_{\Gamma, \Z/2\Z}]$.

Kubjas \cite{Kubjas} and Donten-Bury \cite{DBM} showed that Hilbert functions no longer agree for various other finite abelian groups $G$, so the existence of a common flat deformation cannot hold for the phylogenetic group-based models in general. However, Kubjas and Manon \cite{kubjas2014conformal} have shown that a generalization of the relationship to the WZW model of conformal field theory and the moduli of vector bundles holds for the cyclic groups $\Z/m\Z$. In particular, these group-based models are related to the corresponding moduli spaces for the algebraic group $SL_m(\kk)$.

To sum up, group-based models:
\begin{enumerate}
\item can be regarded as basic combinatorial objects encoding a structure of a finite abelian group,
\item first appeared in phylogenetics,
\item are important also in other fields, such as conformal field theory.
\end{enumerate}

Below we present a table with known facts about generators of ideals for phylogenetic group-based models.

\begin{center}
\begin{tabular}{|p{3,1 cm}|p{1,8 cm}|p{2,3 cm}|p{2,7 cm}|p{4,5 cm}|}
\hline
&\multicolumn{4}{ |c| }{Group-based Models} \\
\hline
Polynomials defining:& $\Z_2$ & $\Z_3$ & $\Z_2\times\Z_2$ & $G$ \\
\hline
Generators of the ideal & Degree $2$ \cite{SS} & Degree $3$ \cite{jaZ3nowy} & Conjecture \cite[Conjecture 30]{SS}  & Finite \cite{jaEmanuele}, Degree $\leq |G|$ \cite[Conjecture 29]{SS} \\
\hline
Projective scheme & & Degree $3$ \cite{Marysianowyz3} & Degree $4$ \cite{JaJCTA} & Finite \cite{JaJCTA}\\
\hline
Set-theoretically&&&&Finite \cite{DE}
\\
\hline
On a Zariski open subset & && Degree $4$ \cite{JaAdvGeom} & Degree $\leq|G|$ \cite{CFSM, casanellas2015complete}\\
\hline
\end{tabular}
\end{center}

\subsection{Open problems}
Here we present the main open problems concerning group-based models. Everything is stated in purely toric/combinatorial language.
We start from the central conjecture in this context. 

\begin{con}[{\cite[Conjecture 29]{SS}}]
For any finite abelian group $G$, $\phi(G)\leq |G|$.
\end{con}
It seems crucial to first understand the simplest tree $K_{1,3}$.
\begin{con}
For any finite abelian group $G$, $\phi(G,3)\leq |G|$.
\end{con}
The results of \cite{jaEmanuele} imply that for finite abelian group $G$ the function $\phi(G,\cdot)$ is eventually constant.
The ensuing results would be a desired strengthening.
\begin{con}[{\cite[Conjecture 9.3]{JaJCTA}}]\label{con:constant}
We have $\phi(G,n+1)=\max(2,\phi(G,n))$.
\end{con}
We are grateful to Seth Sullivant for noticing that this is equivalent to $\phi(G,\cdot)$ being constant, apart from the case when $G=\Z_2$ and $n=3$, when the associated variety is the whole projective space.
\begin{con}
For any finite (not necessarily abelian) group $G$, $\phi(G)$ is finite. 
\end{con}
Conjecture \ref{con:constant} also implies the following.
\begin{con}[{\cite[Conjecture 30]{SS}}]
The phylogenetic complexity of $G=\Z_2 \times \Z_2$ is $4$.
\end{con}

\section{Maps of Toric Varieties and Cox Rings}\label{sec:JB}
Throughout the review we did not mention many important topic, among those Cox rings. In this section we very briefly present the general construction and relations to morphisms.  

Local coordinate rings are not always very convenient to work with, especially,
  when we want to investigate the \emph{global} properties of the variety.
Consider the \emph{projective space} $\p^n$ of dimension $n$.
It is glued out of $n+1$ affine spaces of dimension $n$,
  so to obtain the  description of (for example) a coherent sheaf on the projective space one needs
  the information about $n+1$ modules over polynomial rings,
  and a care must be taken to \emph{glue} the modules accordingly.
Instead, one may view the projective space globally:
\[
  \p^n = (\a^{n+1}\setminus \set{0}) / \C^*.
\]

There are three essential ingredients in this \emph{global description}.
Firstly, there is $\a^{n+1}$, an affine space.
Secondly, we remove a relatively small subset of the affine space, in this case just one point $\{0\}$.
Thirdly, we divide by an action of an algebraic group $\C^*$, the multiplicative group of the base field $\C$.
The \emph{homogeneous coordinate ring} of the projective space incorporates all these three ingredients.
We just take the polynomial coordinate ring of $\a^{n+1}$;
  all objects (for example modules or ideals) that are supported in $\{0\}$ are \emph{irrelevant},
     and, in particular, if two object differ only at $\{0\}$,
     then they correspond to the same object on the projective space;
  all objects must be invariant with respect to the group action, in other words \emph{homogeneous}.
%

The Cox rings have been first introduced for toric varieties \cite{cox_homogeneous},
   and then generalised to normal varieties
   with finitely generated divisor class group $Cl(X)$:
 
\[
   S[X] :=  \bigoplus_{[D] \in Cl(X)} H^0(\mathcal{O}_X(D)).
\]
 
A careful choice of the representatives $D$ in each element of $Cl(X)$ must be made in order to obtain a well defined ring structure on $S[X]$.
Varieties, for which the Cox ring is finitely generated are called \emph{Mori Dream Spaces (MDS)} \cite{hu_keel_MDS_and_GIT}.
The same varieties arise naturally in Mori theory and Minimal Model Program, as particularly elegant examples illustrating the theory.
$S[X]$ is always graded by $Cl(X)$.
The main point is that for a MDS there exists a codimension at least $2$ subvariety $Z$ of $\Spec S[X]$, such that:
\[
  X = (\Spec S[X] \setminus Z)/ G_X
\]
where $G_X = \Hom(Cl(X), \C^*)$ is the group acting on $\Spec S[X]$, corresponding to the grading by $Cl(X)$.
This naturally corresponds to the three ingredients of the homogeneous coordinate ring of the projective space.
Analogously to the projective case, many global objects on $X$ can be expressed in terms of the Cox ring and vice versa,
   taking in account the homogeneity and relevance.
The Cox ring (as well as its grading and the irrelevant ideal $B = I(Z)$ defining $Z$) are defined intrinsically,
   so they do not depend on any embedding or any other choices.
Thus it is very convenient to study \emph{global and intrinsic properties} of $X$.

Affine and projective spaces and normal toric varieties are Mori Dream Spaces.
In these cases, the Cox ring is always a polynomial ring, but the grading vary.
In fact, the property that $S[X]$ is a polynomial ring characterises toric varieties,
  see \cite{kedzierski_wisniewski_Jaczewski_theorem_revisited}
  for a recent treatment of this characterisation.

By definition, an algebraic morphism of two affine varieties $\varphi\colon X\to Y $
  is a geometric interpretation of an algebra morphism $\varphi^* \colon B \rightarrow A$ of their affine coordinate rings.  
Here $X = \Spec A$ and $Y = \Spec B$.
If $X = \p^m$ and $Y = \p^n$ instead,
  and $A \simeq \C[\fromto{x_0}{x_m}]$ and $B \simeq \C[\fromto{y_0}{y_n}]$ are their homogeneous coordinate rings,
  then any algebraic morphism $\varphi\colon \p^m \to \p^n$ is determined a morphism $B \rightarrow A$
  satisfying the usual homogeneity and base point freeness conditions.
Rational maps between affine varieties or projective spaces have similar interpretations in terms of the fields of
  fractions of coordinate rings.


\begin{thm}[{\cite{brown_jabu_maps_of_toric_varieties}, \cite{jabu_kedzierski_maps_of_MDS_1}}]
  Suppose $X$ and $Y$ are Mori Dream Spaces, and $\varphi\colon X \dashrightarrow Y$ is a rational map.
  Then there exists a description of $\varphi$ in terms of Cox coordinates, that is a multi-valued map
  \[
     \Phi \colon \Spec S[X] \multito \Spec S[Y]
  \]
  such that for all points $x \in X$ and  $\xi$ such that $\pi_X(\xi)=x$ and $\varphi$ is regular at $x$,
          the composition $\pi_Y(\Phi(\xi))$ is a single point $ \varphi(x) \in Y$.
\end{thm}
 
The notion of \emph{multi-valued map} is modeled on the case of projective space, 
  but may involve roots of homogeneous functions if the target is singular.
Just as in the case of projective space, the map must satisfy homogeneity, and relevance condition.
The theorem is effective in the sense, that the proof shows how to construct the description.

Similar statement for regular maps between $\q$-factorial Mori Dream Spa\-ces
was obtained by Andreas Hochenegger and Elena Martinengo \cite{hochenegger_martinengo_maps_of_MDS}.
Their approach is to use the language of Mori Dream stacks \cite{hochenegger_martinengo_MD_stacks}.
They use the technique of root constructions, which is parallel to the multi-valued maps.

\section{Examples}\label{sec:ex}
Toric varieties provide very fruitful examples. This section is motivated by questions of Sijong Kwak: what happens to the depth under inner projection of a projective variety $X$? Here, inner projection means a projection from a point $x\in X$. We will denote the (closure of) the image of the projection by $X_x$.
We start by recalling the following general result.
\tw[\cite{Kwak} Theorem 4.1]\label{thm:Kwak}
If $X$ is defined by quadrics and $x$ is a smooth point of $X$ then depths of $X$ and $X_x$ are equal.
\ktw
The following observation was pointed out by Greg Blekherman.
\ex 
In general, the depth may go up under projections from general (in particular, smooth) points. Indeed, if we consider any non aCM variety $X$ we may project it, until it becomes a hypersurfece. In particular, it becomes a complete intersection, hence aCM, hence of maximal depth.
\kex 
Before we pass to constructing toric examples we note that inner projections in toric geometry were investigated for many years. The seminal work of Bruns and Gubeladze \cite{Bruns1, Bruns2} lead to many interesting examples, disproving important conjectures on characterizations of normal polytopes. From a combinatorial point of view projecting from a torus invariant point corresponding to a vertex $v$ of a lattice polytope $P$ corresponds to considering a toric variety given by lattice points in $P$ distinct from $P$. The study when such polytopes remain normal, i.e.~when the projected variety is projectively normal, were crucial in \cite{Bruns1}. Further examples of projectively normal toric varieties that do \emph{not} come from projections of projectively normal toric varieties were found in \cite{BGM}.
Projective normality of toric varieties is related to depth as follows.
\tw[Hochster \cite{HochsterRingsOf}]
A projectively normal toric variety is aCM.
\ktw
\ex 
Consider a toric hypersurface $X$ corresponding to lattice points $(0,0,0)$, $(0,1,0)$, $(0,0,1)$, $(3,1,1)$, $(4,1,1)$.
The Macaulay2 code below verifies that it is aCM and not normal.
\begin{verbatim}
loadPackage "Depth"
loadPackage "Normaliz"
loadPackage("MonomialAlgebras",Configuration=>{"Use4ti2"=>true})
L={{1,0,0,0},{1,0,1,0},{1,0,0,1},{1,4,1,1},{1,3,1,1}}
R=QQ[a_1..a_5,Degrees=>L]
J1=binomialIdeal R
depth (R/J1)==(dim J1)
isNormal (R/J1)
\end{verbatim}
The reason is that the singular locus is of codimension one.
\kex
\ex[The depth may go down under projection from a generic (in particular smooth) point]
We start with a normal (aCM) projective toric variety defined below.
\begin{verbatim}
L={{1,0,0,0},{1,0,1,0},{1,0,0,1},{1,4,1,1},{1,3,1,1},{1,2,1,1},{1,-1,0,0}}
R=QQ[a_1..a_7,Degrees=>L]
J2=binomialIdeal R
isNormal (R/J2)
depth (R/J2)==(dim J2)
\end{verbatim}
We now project from the point $(1,\dots,1)$ that belongs to the dense torus orbit.
\begin{verbatim}
JJ=sub(J2,{a_2=>a_2+a_1,a_3=>a_3+a_1,a_4=>a_4+a_1,
a_5=>a_5+a_1,a_6=>a_6+a_1,a_7=>a_7+a_1})
JS=eliminate(JJ,a_1);
W=QQ[a_1..a_7]
M=sub(JS,W)+ideal(a_1)
depth (W/M)==(dim M)
\end{verbatim}
\kex
In a similar way one can construct examples projecting from singular points. 
\uwa\label{rem:projnonnorm}
It is \emph{not} possible to project a nonprojectively normal toric variety from a (torus invariant) smooth point and obtain a projectively normal variety (as union of normal polytopes is normal). 
\kuwa
We note that the discussion on projections nicely ties with the conjectures of Bogvad and Oda.
\hip\label{hip:BO}
A smooth polytope is normal. The associated toric variety is defined by quadrics.
\khip
\ob
Conjecture \ref{hip:BO} implies that for a smooth polytope $P$ any projection from a torus invariant point remains projectively normal. Further, if we know for a smooth polytope $P$ that there exists a projection from a torus invariant point, that is (projectively) normal, then $P$ is normal.
\kob
\begin{proof}
The second statement follows by Remark \ref{hip:BO}. The first statement is based on \cite[Theorem 5.1]{Bruns:quest} and \cite[Section 11]{jaDissert}. Let $Q$ be the convex hull of lattice points in $P$ distinct from a vertex $v$.
Let $q\in kQ$. We know that $q=\sum_{i=1}^k p_i$ for lattice points $p_i\in P$. Let $v_1,\dots, v_n$ be the first lattice points on edges of $P$ adjacent to $v$. The only problem is, if some $p_i=v$. But then there must exist $p_j\neq v,v_1,\dots,v_n$. Note that as $v-v_i$ is a lattice basis we may always find $m\in \z_+$ such that $mv+p_j=\sum v_i$, where the sum is over any indices (with possible repetitions) $1\leq i\leq n$. The only way the above relation is generated by quadrics is, if $v+p_j=a+b$ for some lattice points $a,b\in P$. Thus as long as in the decomposition of $q$ the vertex $v$ appears we may change the decomposition in such a way that its multiplicity goes down. This proves normality of $Q$. 
\end{proof}
The previous proposition ties with the general Theorem \ref{thm:Kwak}. Indeed, if $P$ is defined by quadrics and normal, then it is aCM and we know that the projection is also aCM. 

\bibliographystyle{amsalpha}
\bibliography{Xbib}

\newcommand{\etalchar}[1]{$^{#1}$}
\def\polhk#1{\setbox0=\hbox{#1}%
  {\ooalign{\hidewidth\lower1.5ex\hbox{`}\hidewidth\crcr\unhbox0}}}\def\dbar{\leavevmode\hbox
  to 0pt{\hskip.2ex\accent"16\hss}d}
\providecommand{\bysame}{\leavevmode\hbox to3em{\hrulefill}\thinspace}
\providecommand{\MR}{\relax\ifhmode\unskip\space\fi MR }
\providecommand{\MRhref}[2]{%
  \href{http://www.ams.org/mathscinet-getitem?mr=#1}{#2}
}
\providecommand{\href}[2]{#2}
\begin{thebibliography}{CFSM15b}

\bibitem[AIP{\etalchar{+}}11]{altmann2011geometry}
Klaus Altmann, Nathan~O Ilten, Lars Petersen, Hendrik S{\"u}{\ss}, and Robert
  Vollmert, \emph{The geometry of t-varieties}, Contributions to Algebraic
  Geometry, Impanga Lect Notes (2011).

\bibitem[AM69]{AtiMac}
Michael~Francis Atiyah and Ian~Grant Macdonald, \emph{Introduction to
  commutative algebra}, vol.~2, Addison-Wesley Reading, 1969.

\bibitem[BB13]{brown_jabu_maps_of_toric_varieties}
Gavin Brown and Jaros\l{}aw Buczy\'nski, \emph{Maps of toric varieties in {C}ox
  coordinates}, Fund. Math. \textbf{222} (2013), 213--267.

\bibitem[BBD{\etalchar{+}}16]{banos2016phylogenetic}
Hector Ba{\~n}os, Nathaniel Bushek, Ruth Davidson, Elizabeth Gross, Pamela~E
  Harris, Robert Krone, Colby Long, Allen Stewart, and Robert Walker,
  \emph{Phylogenetic trees}, arXiv preprint arXiv:1611.05805 (2016).

\bibitem[BG99]{Bruns1}
Winfried Bruns and Joseph Gubeladze, \emph{Normality and covering properties of
  affine semigroups}, J. Reine Angew. Math. \textbf{510} (1999), 161--178.
  \MR{1696094}

\bibitem[BG09]{BG}
\bysame, \emph{Polytopes, rings, and k-theory}, Springer Science \& Business
  Media, 2009.

\bibitem[BGM16]{BGM}
Winfried Bruns, Joseph Gubeladze, and Mateusz Micha{\l}ek, \emph{Quantum jumps
  of normal polytopes}, Discrete Comput. Geom. \textbf{56} (2016), no.~1,
  181--215. \MR{3509036}

\bibitem[BGT97]{Bruns2}
Winfried Bruns, Joseph Gubeladze, and Ng{\^o}~Vi{\^e}t Trung, \emph{Normal
  polytopes, triangulations, and {K}oszul algebras}, J. Reine Angew. Math.
  \textbf{485} (1997), 123--160. \MR{1442191}

\bibitem[BH09]{borisov2009conjecture}
Lev Borisov and Zheng Hua, \emph{On the conjecture of king for smooth toric
  deligne--mumford stacks}, Advances in Mathematics \textbf{221} (2009), no.~1,
  277--301.

\bibitem[BIS]{Normaliz}
Winfried Bruns, Bogdan Ichim, and Christof S\"{o}ger, \emph{Normaliz},
  http://www.mathematik.uni-osnabrueck.de/normaliz.

\bibitem[BK16]{jabu_kedzierski_maps_of_MDS_1}
Jaros{\l}aw Bucz{y\'n}ski and Oskar K{\k e}dzierski, \emph{Maps of mori dream
  spaces in cox coordinates. part i: existence of descriptions}, arXiv:
  1605.06828, 2016.

\bibitem[Bla08]{Blasiak}
Jonah Blasiak, \emph{The toric ideal of a graphic matroid is generated by
  quadrics}, Combinatorica \textbf{28} (2008), no.~3, 283--297.

\bibitem[Blu01]{Bl01}
Stefan Blum, \emph{Base-sortable matroids and koszulness of semigroup rings},
  European Journal of Combinatorics \textbf{22} (2001), no.~7, 937--951.

\bibitem[BM86]{CutPolytope}
Francisco Barahona and Ali~Ridha Mahjoub, \emph{On the cut polytope},
  Mathematical programming \textbf{36} (1986), no.~2, 157--173.

\bibitem[Bon13]{Bonin}
Joseph~E Bonin, \emph{Basis-exchange properties of sparse paving matroids},
  Advances in Applied Mathematics \textbf{50} (2013), no.~1, 6--15.

\bibitem[Bru13]{Bruns:quest}
Winfried Bruns, \emph{The quest for counterexamples in toric geometry},
  Commutative algebra and algebraic geometry ({CAAG}-2010), Ramanujan Math.
  Soc. Lect. Notes Ser., vol.~17, Ramanujan Math. Soc., Mysore, 2013,
  pp.~45--61. \MR{3155951}

\bibitem[Buc12]{MR2892983}
Weronika Buczy{\'n}ska, \emph{Phylogenetic toric varieties on graphs}, J.
  Algebraic Combin. \textbf{35} (2012), no.~3, 421--460. \MR{2892983}

\bibitem[BW07]{BW}
Weronika Buczy\'nska and Jaros{\l}aw~A. Wi\'{s}niewski, \emph{On geometry of
  binary symmetric models of phylogenetic trees}, J. Eur. Math. Soc.
  \textbf{9(3)} (2007), 609--635.

\bibitem[CFSM15a]{casanellas2015complete}
Marta Casanellas, Jes{\'u}s Fern{\'a}ndez-S{\'a}nchez, and Mateusz Micha{\l}ek,
  \emph{Complete intersection for equivariant models}, Preprint
  arXiv:1512.07174 (2015).

\bibitem[CFSM15b]{CFSM}
\bysame, \emph{Low degree equations for phylogenetic group-based models},
  Collect. Math. \textbf{66} (2015), no.~2, 203--225.

\bibitem[CLS11]{Cox}
David~A Cox, John~B Little, and Henry~K Schenck, \emph{Toric varieties},
  American Mathematical Soc., 2011.

\bibitem[Con07]{Conca}
Aldo Conca, \emph{Linear spaces, transversal polymatroids and asl domains},
  Journal of Algebraic Combinatorics \textbf{25} (2007), no.~1, 25--41.

\bibitem[Cox95]{cox_homogeneous}
David~A. Cox, \emph{The homogeneous coordinate ring of a toric variety}, J.
  Algebraic Geom. \textbf{4} (1995), no.~1, 17--50. \MR{MR1299003 (95i:14046)}

\bibitem[DB16]{Marysianowyz3}
Maria Donten-Bury, \emph{Phylogenetic invariants for $\mathbb{Z}_3$
  scheme-theoretically}, Ann. Comb. \textbf{20} (2016), no.~3, 549--568.

\bibitem[DBM12]{DBM}
Maria Donten-Bury and Mateusz Micha{\l}ek, \emph{Phylogenetic invariants for
  group-based models}, J. Algebr. Stat. \textbf{3} (2012), no.~1, 44--63.

\bibitem[DE15]{DE}
Jan Draisma and Rob~H. Eggermont, \emph{Finiteness results for {A}belian tree
  models}, J. Eur. Math. Soc. (JEMS) \textbf{17} (2015), no.~4, 711--738.

\bibitem[Eng11]{Engstrom}
Alexander Engstr{\"o}m, \emph{Cut ideals of k4-minor free graphs are generated
  by quadrics}, Michigan Math. J. \textbf{60} (2011), no.~3, 705--714.

\bibitem[ERSS05]{4aut}
Nicholas Eriksson, Kristian Ranestad, Bernd Sturmfels, and Seth Sullivant,
  \emph{Phylogenetic algebraic geometry}, Projective varieties with unexpected
  properties (2005), 237--255.

\bibitem[Fed83]{Fedder}
Richard Fedder, \emph{F-purity and rational singularity}, Transactions of the
  American Mathematical Society \textbf{278} (1983), no.~2, 461--480.

\bibitem[Ful93]{Ful}
William Fulton, \emph{Introduction to {Toric} {Varieties}}, Annals of
  Mathematics Studies, vol. 131, Princeton University Press, 1993.

\bibitem[GGMS87]{GGMS}
Israel~M Gelfand, R~Mark Goresky, Robert~D MacPherson, and Vera~V Serganova,
  \emph{Combinatorial geometries, convex polyhedra, and schubert cells},
  Advances in Mathematics \textbf{63} (1987), no.~3, 301--316.

\bibitem[GJ00]{polymake}
Ewgenij Gawrilow and Michael Joswig, \emph{Polymake: a {Framework} for
  {Analyzing} {Convex} {Polytopes}}, Polytopes --- Combinatorics and
  Computation (Gil Kalai and G\"unter~M. Ziegler, eds.), Birkh\"auser, 2000,
  pp.~43--74.

\bibitem[GR12]{GiRe12}
Dion Gijswijt and Guus Regts, \emph{Polyhedra with the integer caratheodory
  property}, Journal of Combinatorial Theory, Series B \textbf{102} (2012),
  no.~1, 62--70.

\bibitem[GS]{M2}
Daniel~R. Grayson and Michael~E. Stillman, \emph{Macaulay2, a software system
  for research in algebraic geometry}, Available at
  http://www.math.uiuc.edu/Macaulay2/.

\bibitem[Har13]{Hart}
Robin Hartshorne, \emph{Algebraic geometry}, vol.~52, Springer Science \&
  Business Media, 2013.

\bibitem[HH02]{HerzogHibi}
J{\"u}rgen Herzog and Takayuki Hibi, \emph{Discrete polymatroids}, Journal of
  Algebraic Combinatorics \textbf{16} (2002), no.~3, 239--268.

\bibitem[HHV{\etalchar{+}}05]{HHV}
J{\"u}rgen Herzog, Takayuki Hibi, Marius Vladoiu, et~al., \emph{Ideals of fiber
  type and polymatroids}, Osaka Journal of Mathematics \textbf{42} (2005),
  no.~4, 807--829.

\bibitem[HK00]{hu_keel_MDS_and_GIT}
Yi~Hu and Sean Keel, \emph{Mori dream spaces and {GIT}}, Michigan Math. J.
  \textbf{48} (2000), 331--348, Dedicated to William Fulton on the occasion of
  his 60th birthday. \MR{1786494 (2001i:14059)}

\bibitem[HK12]{Kwak}
Kangjin Han and Sijong Kwak, \emph{Analysis on some infinite modules, inner
  projection, and applications}, Transactions of the American Mathematical
  Society \textbf{364} (2012), no.~11, 5791--5812.

\bibitem[HM15]{hochenegger_martinengo_MD_stacks}
Andreas Hochenegger and Elena Martinengo, \emph{Mori dream stacks}, Math. Z.
  \textbf{280} (2015), no.~3-4, 1185--1202. \MR{3369373}

\bibitem[HM16]{hochenegger_martinengo_maps_of_MDS}
\bysame, \emph{Maps of {M}ori {D}ream {S}paces}, arXiv:1605.06789, 2016.

\bibitem[Hoc72]{HochsterRingsOf}
M.~Hochster, \emph{Rings of invariants of tori, {C}ohen-{M}acaulay rings
  generated by monomials, and polytopes}, Ann. of Math. (2) \textbf{96} (1972),
  318--337. \MR{0304376}

\bibitem[HR76]{HochsterRoberts}
Melvin Hochster and Joel~L Roberts, \emph{The purity of the frobenius and local
  cohomology}, Advances in Mathematics \textbf{21} (1976), no.~2, 117--172.

\bibitem[Kas10]{Ka10}
Kenji Kashiwabara, \emph{The toric ideal of a matroid of rank 3 is generated by
  quadrics}, Electron. J. Combin \textbf{17} (2010), no.~1.

\bibitem[Kim81]{3Kimura}
Motoo Kimura, \emph{Estimation of evolutionary distances between homologous
  nucleotide sequences}, Proceedings of the National Academy of Sciences
  \textbf{78} (1981), no.~1, 454--458.

\bibitem[KKMSD]{kempt339toroidal}
G~Kempf, F~Knudson, D~Mumford, and B~Saint-Donat, \emph{Toroidal embedding 1},
  Lecture Notes in Math \textbf{339}.

\bibitem[Kly90]{klyachko1990equivariant}
Alexander~A Klyachko, \emph{Equivariant bundles on toral varieties},
  Mathematics of the USSR-Izvestiya \textbf{35} (1990), no.~2, 337.

\bibitem[KM14]{kubjas2014conformal}
Kaie Kubjas and Christopher Manon, \emph{Conformal blocks,
  {B}erenstein--{Z}elevinsky triangles, and group-based models}, J. Algebr.
  Comb. \textbf{40} (2014), no.~3, 861--886.

\bibitem[Kub12]{Kubjas}
Kaie Kubjas, \emph{Hilbert polynomial of the kimura 3-parameter model.},
  Journal of Algebraic Statistics \textbf{3} (2012), no.~1.

\bibitem[KW15]{kedzierski_wisniewski_Jaczewski_theorem_revisited}
Oskar K{\polhk{e}}dzierski and Jaros{\l}aw~A. Wi{\'s}niewski,
  \emph{Differentials of {C}ox rings: {J}aczewski's theorem revisited}, J.
  Math. Soc. Japan \textbf{67} (2015), no.~2, 595--608. \MR{3340188}

\bibitem[Las]{La17}
Micha{\l} Laso{\'n}, \emph{On gr \"{o}bner bases of the toric ideal of a
  matroid}, in preparation.

\bibitem[Las15]{La15}
\bysame, \emph{List coloring of matroids and base exchange properties},
  European Journal of Combinatorics \textbf{49} (2015), 265--268.

\bibitem[Las16]{La16}
\bysame, \emph{Degree bounds for the toric ideal of a matroid}, arXiv preprint
  arXiv:1601.08199 (2016).

\bibitem[LM11]{lason2011full}
Micha{\l} Laso{\'n} and Mateusz Micha{\l}ek, \emph{On the full, strongly
  exceptional collections on toric varieties with picard number three},
  Collectanea mathematica \textbf{62} (2011), no.~3, 275--296.

\bibitem[LM14]{jaMichal}
\bysame, \emph{On the toric ideal of a matroid}, Advances in Mathematics
  \textbf{259} (2014), 1--12.

\bibitem[Man12]{M}
Christopher Manon, \emph{Coordinate rings for the moduli stack of
  {$SL_2(\Bbb{C})$} quasi-parabolic principal bundles on a curve and toric
  fiber products}, J. Algebra \textbf{365} (2012), 163--183. \MR{2928457}

\bibitem[Mic11]{JaJalg}
Mateusz Micha{\l}ek, \emph{Geometry of phylogenetic group-based models}, J.
  Algebra \textbf{339} (2011), no.~1, 339--356.

\bibitem[Mic13]{JaJCTA}
\bysame, \emph{Constructive degree bounds for group-based models}, J. Combin.
  Theory, Series A \textbf{120} (2013), no.~7, 1672--1694.

\bibitem[Mic14]{JaAdvGeom}
\bysame, \emph{Toric geometry of the 3-{K}imura model for any tree}, Adv. Geom.
  \textbf{14} (2014), no.~1, 11--30.

\bibitem[Mic15]{jaDissert}
\bysame, \emph{Toric varieties in phylogenetics}, Diss. Math. (2015), no.~511,
  3--86.

\bibitem[Mic17]{jaZ3nowy}
Mateusz Micha{\l}ek, \emph{Finite phylogenetic complexity of $\mathbb {Z}_p$
  and invariants for $\mathbb {Z}_3$}, European J. Combin. \textbf{59} (2017),
  169 -- 186.

\bibitem[MOZ14]{MOZ}
Mateusz Micha{\l}ek, Luke Oeding, and Piotr Zwiernik, \emph{Secant cumulants
  and toric geometry}, International Mathematics Research Notices (2014),
  rnu056.

\bibitem[MS15]{maclagan2015introduction}
Diane Maclagan and Bernd Sturmfels, \emph{Introduction to tropical geometry},
  vol. 161, American Mathematical Soc., 2015.

\bibitem[MV17]{jaEmanuele}
Mateusz Micha{\l}ek and Emanuele Ventura, \emph{Finite phylogenetic complexity
  and combinatorics of tables}, Algebra and Number Theory \textbf{11} (2017),
  no.~1.

\bibitem[Oda88]{Oda}
Tadao Oda, \emph{Convex bodies and algebraic geometry}, Ergebnisse der
  Mathematik und ihrer Grenzgebiete (3) [Results in Mathematics and Related
  Areas (3)], vol.~15, Springer-Verlag, Berlin, 1988, An introduction to the
  theory of toric varieties, Translated from the Japanese. \MR{MR922894
  (88m:14038)}

\bibitem[OH99]{ohsugi1999toric}
Hidefumi Ohsugi and Takayuki Hibi, \emph{Toric ideals generated by quadratic
  binomials}, Journal of Algebra \textbf{218} (1999), no.~2, 509--527.

\bibitem[Ohs10]{OhsugiCuts}
Hidefumi Ohsugi, \emph{Normality of cut polytopes of graphs is a minor closed
  property}, Discrete Mathematics \textbf{310} (2010), no.~6, 1160--1166.

\bibitem[Oxl06]{Ox}
James~G Oxley, \emph{Matroid theory}, vol.~3, Oxford University Press, USA,
  2006.

\bibitem[Pay09]{payne2009toric}
Sam Payne, \emph{Toric vector bundles, branched covers of fans, and the
  resolution property}, Journal of Algebraic Geometry \textbf{18} (2009),
  no.~1, 1--36.

\bibitem[Sch11]{Sc11}
Jay Schweig, \emph{Toric ideals of lattice path matroids and polymatroids},
  Journal of Pure and Applied Algebra \textbf{215} (2011), no.~11, 2660--2665.

\bibitem[Sey81]{SeymourMatroids}
PD~Seymour, \emph{Matroids and multicommodity flows}, European Journal of
  Combinatorics \textbf{2} (1981), no.~3, 257--290.

\bibitem[Shi16]{shibata2016toric}
Kazuki Shibata, \emph{Toric ideals of series and parallel connections of
  matroids}, Journal of Algebra and Its Applications \textbf{15} (2016),
  no.~06, 1650106.

\bibitem[SS05]{SS}
Bernd Sturmfels and Seth Sullivant, \emph{Toric ideals of phylogenetic
  invariants}, J. Comput. Biology \textbf{12} (2005), 204--228.

\bibitem[SS08]{Cutsandsplits}
Bernd Sturmfels and Seth Sullivant, \emph{Toric geometry of cuts and splits},
  The Michigan Mathematical Journal \textbf{57} (2008), 689--709.

\bibitem[Stu96]{Stks}
Bernd Sturmfels, \emph{Gr\"obner bases and convex polytopes}, University
  Lecture Series, vol.~8, American Mathematical Society, 1996.

\bibitem[Sul07]{Sethtfp}
Seth Sullivant, \emph{Toric fiber products}, J. Algebra \textbf{316} (2007),
  no.~2, 560 -- 577.

\bibitem[Sum74]{sumihiro1974equivariant}
Hideyasu Sumihiro, \emph{Equivariant completion}, Journal of Mathematics of
  Kyoto University \textbf{14} (1974), no.~1, 1--28.

\bibitem[SX10]{SX}
Bernd Sturmfels and Zhiqiang Xu, \emph{Sagbi bases of {C}ox-{N}agata rings}, J.
  Eur. Math. Soc. (JEMS) \textbf{12} (2010), no.~2, 429--459. \MR{2608947}

\bibitem[tt]{4ti2}
4ti2 team, \emph{4ti2---a software package for algebraic, geometric and
  combinatorial problems on linear spaces}, www.4ti2.de.

\bibitem[Whi77]{white1977basis}
Neil~L White, \emph{The basis monomial ring of a matroid}, Advances in
  Mathematics \textbf{24} (1977), no.~2, 292--297.

\bibitem[Whi80]{White}
\bysame, \emph{A unique exchange property for bases}, Linear Algebra and its
  applications \textbf{31} (1980), 81--91.

\end{thebibliography}
\noindent Mateusz Micha{\l}ek,
Polish Academy of Sciences, Warsaw, Poland,\\
Research Institute for Matematical Sciences, Kyoto, Japan\\
{\tt wajcha2@poczta.onet.pl}
\end{document}